\documentclass[12pt]{amsart}

\usepackage{pb-diagram}

\usepackage{xcolor}

\usepackage{amsfonts}
\usepackage{amsmath}
\usepackage{amssymb}

\allowdisplaybreaks

\newcommand{\cS}{\mathcal{S}}

\newcommand{\mref}[1]{(\ref{#1})}

\newtheorem{theorem}{Theorem}[section]

\renewcommand{\thefootnote}{\fnsymbol{footnote}}
\newcommand{\starttext}{ \setcounter{footnote}{0}
\renewcommand{\thefootnote}{\arabic{footnote}}}

\newcommand{\beq}{\begin{equation}}
\newcommand{\bea}{\begin{eqnarray}}
\newcommand{\eea}{\end{eqnarray}} \newcommand{\ee}{\end{equation}}
 \newcommand{\<}{\langle}
\renewcommand{\>}{\rangle}
\def\ba{\begin{eqnarray}}
\def\ea{\end{eqnarray}}

\def\o{\omega}

\def\ti\tilde

\def\u{\underline}

\def\pl{\partial}
\def\na{\nabla}

\def\p{\prod}

\def\ddb{\partial\bar\partial}

\def\Im{{\rm Im}}

\def\ti{\tilde}

\begin{document}

\starttext \baselineskip=18pt \setcounter{footnote}{0}

\starttext \baselineskip=18pt \setcounter{footnote}{0}

\title[Geometric invariants and the Monge-Amp\`ere equation]{\Large GEOMETRIC INVARIANTS AND THE MONGE-AMP\`ERE EQUATION IN K\"AHLER GEOMETRY}

\author{Bin Guo  and Duong H. Phong }

\thanks{Work supported in part by the National Science Foundation under grants DMS-22-03273, DMS-25-05241,  and DMS-23-03508, and the collaboration grant 946730 from Simons Foundation.}

\address{ Department of Mathematics \& Computer Science, Rutgers University, Newark, NJ 07102}

\email{bguo@rutgers.edu}

\address{Department of Mathematics, Columbia University, New York, NY 10027}

\email{phong@math.columbia.edu}

\begin{abstract}

This is a contribution to the special issue of Surveys in Differential Geometry celebrating the 75th birthday of Shing-Tung Yau. The bulk of the paper is devoted to a survey of some new geometric inequalities and estimates for the Monge-Amp\`ere equation, obtained by the authors in the last few years in joint work with F. Tong, J. Song, and J. Sturm. These all depend 
in an essential way on Yau's solution of the Calabi conjecture, which is itself nearing its own 50th birthday. The opportunity is also taken to survey briefly many current directions in complex geometry, which he more recently pioneered.

\end{abstract}

{\footnotesize}

\maketitle


\baselineskip=15pt
\setcounter{equation}{0}
\setcounter{footnote}{0}

\section{\bf INTRODUCTION}
\setcounter{equation}{0}

It is a great honor and a great pleasure to contribute this paper to the issue of Surveys in Differential Geometry celebrating the 75th birthday of Shing-Tung Yau. As we are also approaching the 50th anniversary of his landmark solution of the Calabi conjecture, it may be appropriate to concentrate the bulk of the present survey on some recent applications of this solution. More specifically, many new inequalities for the diameter, volume density, Green's function, Sobolev inequality, and Gromov convergence theorems for K\"ahler manifolds have been obtained using the method of auxiliary Monge-Amp\`ere equations, which relies in an essential way on Yau's solution of the Calabi conjecture. What makes these new inequalities special is their uniformity with respect to large classes of K\"ahler metrics, and the perhaps surprising fact that they
do not require any lower bound assumption on the Ricci curvature. As such, they have no known analogue for Riemannian manifolds, and seem to open up many possibilities which may well alter the landscape of K\"ahler geometry.

\bigskip
Ricci-flat K\"ahler metrics, whose existence had been conjectured by Calabi and proved by Yau \cite{Y} to exist on K\"ahler manifolds with vanishing first Chern class, are now known as Calabi-Yau metrics and their underlying space as Calabi-Yau spaces. They can be viewed as the higher-dimensional analogue of elliptic curves, and just as elliptic curves had revealed themselves to be at the crossroads of practically every branch of mathematics, so have Calabi-Yau spaces. Calabi-Yau spaces are arguably even more fundamental, as they have also turned out to be essential in the understanding of string theories, presently still the only viable candidate for a unified theory of all physical interactions \cite{CHSW}. Thus any survey of the ramifications of Yau's solution of the Calabi conjecture is necessarily woefully inadequate. Our choice of developments closely related to our areas of expertise is only a reflection of our own particular knowledge and limitation.
Nevertheless, the other developments are also extremely exciting and will undoubtedly prove to be a fertile ground for research in the foreseeable future. With this in mind,
we thought that it would be beneficial to young students and researchers to include in this survey an introduction, even if it is not much more than a guide to the literature, to some of the developments in complex geometry and partial differential equations which Yau had himself pioneered on topics closely related to Calabi-Yau manifolds and canonical metrics. Such an incomplete and sketchy introduction may still be on some value, given the vast expanse of these recent developments, and the accelerating pace of progress in their areas.

\bigskip
The authors would like to thank the editors of this Surveys in Differential Geometry volume for their invitation to contribute to this volume. They would also like to express their gratitude to their many collaborators on these and related issues, notably Freid Tong, Jian Song, Jacob Sturm, Sebastien Picard, Xiangwen Zhang, Teng Fei, and Tristan Collins. The second author would also like to thank Li-Sheng Tseng for explaining to him his joint works with Yau on symplectic cohomology theories.

\bigskip

\section{\bf THE MONGE-AMP\`ERE EQUATION IN K\"AHLER GEOMETRY}
\setcounter{equation}{0}

In this chapter, we recall some basic facts of K\"ahler geometry, the Calabi conjecture, and Yau's Theorem. We then introduce the key new tool leading to the new geometric estimates, namely the method of auxiliary Monge-Amp\`ere equations, and how it applies to Green's functions. As we shall see, this method depends in an essential way on Yau's Theorem. Geometric applications, in particular to a more flexible version of the Gromov convergence theorem and to diameter estimates in the K\"ahler-Ricci flow, are described in the last section.

\subsection{\bf YAU'S THEOREM}

\hfill\break

Let $(X,\o)$ be a compact K\"ahler manifold of dimension $n$, $\o=ig_{\bar kj}dz^j\wedge d\bar z^k$ in local holomorphic coordinates $z^j$, $1\leq j\leq n$. We let $\na$ denote the Chern unitary connection with respect to $\o$,
\bea
\na_{\bar k}V^p=\pl_{\bar k}V^p,
\qquad
\na_jV^p=g^{p\bar m}\pl_j(g_{\bar mq}V^q).
\eea
The curvature tensor of $\na$ is then given by
\bea
[\na_j,\na_{\bar k}]V^p=R_{\bar kj}{}^p{}_qV^q,
\qquad R_{\bar kj}{}^p{}_q=-\pl_{\bar k}(g^{p\bar m}\pl_j g_{\bar mq}).
\eea
The Ricci curvature $R_{\bar kj}$ is defined by
\bea
R_{\bar kj}=R_{\bar kj}{}^p{}_p.
\eea
From the previous formula for the curvature tensor $R_{\bar kj}{}^p{}_q$, it follows immediately that
\bea
R_{\bar kj}=-\pl_{\bar k}(g^{p\bar m}\pl_j g_{\bar mq})
=
-\pl_{\bar k}\pl_j\log \,{\rm det}\,g_{\bar qp}.
\eea
Introducing the Ricci curvature form $Ric(\o)$ by
\bea
Ric(\o)=iR_{\bar kj}dz^j\wedge d\bar z^k
\eea
we can rewrite the previous formula for the Ricci curvature form as
\bea\label{eqn:ric}
Ric(\o)=-i\pl\bar\pl \log\,\o^n.
\eea

\bigskip
Thus $Ric(\o)$ is a closed $(1,1)$-form. If $\tilde\o$ is any other K\"ahler metric on $X$, then
\bea
Ric(\o)-Ric(\tilde\o)=-i\pl\bar\pl\,\log\,{\o^n\over \tilde\o^n}.
\eea
The right hand side is an exact form, since $\o^n/\tilde\o^n$ is a globally defined scalar function. It follows that the Ricci curvature form $Ric(\o)$ depends on $\o$, but its de Rham cohomology class $[Ric(\o)]$ does not, and defines an invariant of the manifold $X$, denoted by $c_1(X)$ and called its first Chern class,
\bea
[Ric(\o)]=c_1(X).
\eea

\bigskip
We can now consider the following K\"ahler version of Einstein's equation in general relativity. Given a form $T=iT_{\bar kj}dz^j\wedge d\bar z^k$, find a K\"ahler metric $\o$
satisfying
\bea
\label{MA0}
Ric(\o)=T.
\eea
In view of the preceding discussion, clearly a necessary condition for the existence of a solution is that the form $T$ be closed, and that its de Rham cohomology class satisfy
\bea
c_1(X)=[T].
\eea
The Calabi conjecture is that this condition is sufficient. More precisely,

\bigskip
{\bf \small Calabi conjecture}: Let $[\o_0]$ be any K\"ahler class on $X$. Then for any given smooth $(1,1)$-form $T$ with $[T]=c_1(X)$, there exists a unique K\"ahler form $\o\in [\o_0]$ satisfying the equation $Ric(\o)=T$.

\bigskip
This was a very bold conjecture, since it implied in particular an abundance of Ricci-flat K\"ahler metrics on K\"ahler manifolds with $c_1(X)=0$, in fact a unique Ricci-flat K\"ahler metric in each K\"ahler class $[\o_0]$. At that time, no Ricci-flat K\"ahler metric was even known except in trivial cases.

\bigskip
The Calabi conjecture can be reduced to the solvability of a Monge-Amp\`ere equation as follows. By the $\pl\bar\pl$-lemma, any K\"ahler metric $\o\in[\o_0]$ can be expressed as
\bea
\o=\o_0+i\pl\bar\pl\varphi
\eea
for a smooth function $\varphi$ which is unique up to an additive constant and is known as the potential of $\o$. Trivially,
\bea
\label{MA1}
(\o_0+i\pl\bar\pl\varphi)^n=e^F\,\o_0^n, \qquad \o_0+i\pl\bar\pl\varphi>0,
\eea
where the relative volume form $F$ is defined by
\bea
F=\log\,{\o^n\over\o_0^n}.
\eea
Conversely, if a relative volume form $F$ is assigned, then the K\"ahler metric $\o$ is determined by its potential $\varphi$, which is itself determined by the Monge-Amp\`ere equation (\ref{MA1}). 

\medskip
Returning to the equation (\ref{MA0}), the basic observation is that the assignment of the Ricci form $Ric(\o)$ is exactly equivalent to the assignment of the relative volume form $F$. Indeed, since $T$ and $Ric(\o_0)$ are both in $c_1(X)$, $T-Ric(\o_0)$ is cohomologically trivial and can be expressed by the $\pl\bar\pl$-lemma as
\bea
T-Ric(\o_0)=-i\pl\bar\pl f
\eea
for a scalar function $f$ which is unique up to an additive constant. But then
\bea
Ric(\o)=T=Ric(\o_0)-i\pl\bar\pl f
\eea
and hence
\bea
-i\pl\bar\pl \,F=Ric(\o)-Ric(\o_0)=-i\pl\bar\pl f
\eea
so that $F=f$ up to an additive constant. This constant is determined by the requirement that
\bea
\label{compatibility}
\int_X \o_0^n=\int_X \o^n=\int_X e^F\o_0^n.
\eea
Thus the Calabi conjecture is equivalent to the solvability of the equation (\ref{MA0}) for an assigned smooth function $F$ satisfying the compatibility condition (\ref{compatibility}). This was proved by Yau \cite{Y}:

\bigskip
{\bf \small Yau's Theorem}: {\it Let $X$ be a compact $n$-dimensional K\"ahler manifold, and let $[\o_0]$ be any given K\"ahler class on $X$. Then 

{\rm (a)} For any smooth strictly positive function $e^F$ on $X$, the Monge-Amp\`ere equation (\ref{MA1}) admits a smooth solution $\varphi$ which is unique up to an additive constant.

{\rm (b)} In particular, for any closed smooth $(1,1)$-form $T$ satisfying $[T]=c_1(X)$, there exists a unique K\"ahler metric $\o\in[\o_0]$ satisfying $Ric(\o)=T$.}

\bigskip
Yau's proof was by the method of continuity. In such a proof, the key step is to obtain suitable {\it a priori estimates} for the equation (\ref{MA1}). More specifically, assume that we have a smooth solution $\varphi$ of the equation (\ref{MA1}), say normalized to satisfy $\int_X\varphi\o_0^n=0$. Then Yau proves:

\medskip
{\it The ($C^0$) estimate}: for any $p>n$, the solution $\varphi$ satisfies
\bea
\|\varphi\|_{C^0}\leq K_{0,p}
\eea
for some constant $K_{0,p}$ which depends only on $X,\o_0$, and $\|e^F\|_{L^p}$.

\medskip
{\it The ($C^2$) estimate}: there exists a constant $K_2$ and a constant $A$, both depending only on $X,\o_0$ and $
\|F\|_{C^2}$ so that
\bea
\Delta_0\varphi\leq K_2\,e^{A\,{\rm osc}\,\varphi}
\eea
where ${\rm osc}\,\varphi={\rm sup}_{z,w\in X}(\varphi(z)-\varphi(w))$, and $\Delta_0\varphi=g_0^{j\bar k}\pl_{\bar k}\pl_j\varphi$ is the Laplacian of $\varphi$ with respect to the metric $\o_0=i(g_0)_{\bar kj}dz^j\wedge d\bar z^k$.

\medskip
{\it The ($C^3$) estimate}: set $S=g_0^{m\bar r}g_0^{s\bar k}g_0^{j\bar t} \na_m^0\pl_{\bar k}\pl_j\varphi\na_{\bar r}^0\pl_s\pl_{\bar t}\varphi$, where $\na^0$ denote the covariant derivatives with respect to the metric $\o_0$. Then 
\bea
S\leq K_3
\eea
for some constant $K_3$ which depends only on $X,\o_0$, $\|F\|_{C^3}$, and ${\rm inf}_X F$.

\bigskip
From these estimates, using the theory of elliptic linear partial differential equations, it is not difficult to establish the following:

\medskip
{\it The ($C^k$) estimate}: assume that $\varphi\in C^5(X)$. Then for any non-negative integer $k\geq 4$, and any $\beta$ with $0<\beta<1$,
there exists a constant $A_{k,\beta}$ depending only on $X,\o_0$, $\|F\|_{C^k}$ and ${\rm inf}_XF$, so that the solution $\varphi$ of the equation (\ref{MA1}) must satisfy
\bea
\|\varphi\|_{C^{k+1,\beta}}\leq A_{k,\beta}.
\eea

\bigskip
\subsection{\bf GEOMETRIC INVARIANTS}

\hfill\break

As discussed in the Introduction, we shall be particularly interested in estimates for the diameter, volume density, Green's function, Sobolev constants, and a Gromov convergence theorem for K\"ahler manifolds. We begin by recalling some basic facts in the general case of a compact Riemannian manifold.

\bigskip
Let $(X,g_{kj})$ be a compact Riemannian manifold. Then its diameter ${\rm diam}_gX$ is defined by
\bea
{\rm diam}_gX={\rm sup}_{x,y\in X}d(x,y)
\eea
where $d(x,y)$ is the distance between $x$ and $y$ with respect to the metric $g_{kj}$.
The Green's function $G_g(x,y)$ is defined by the equations
\bea
\Delta_gG_g(x,\cdot)=-\delta(x,\cdot)+{1\over V_g},
\qquad \int_X G_g(x,\cdot)\sqrt g dx=0,
\eea
where $\sqrt g dx$ is the volume form of $g_{kj}$, $\Delta_g=-{1\over\sqrt g}\pl_j(\sqrt g g^{kj}\pl_k)$ is the Laplacian, and $V_g=\int_X\sqrt g dx$ is the volume of $X$ with respect to the metric $g_{kj}$. The existence and uniqueness of $G_g(x,y)$ follows readily from the theory of linear elliptic partial differential equations (e.g. \cite{ScYa}), together with the fact that it is symmetric $G_g(x,y)=G_g(y,x)$, is smooth off the diagonal $x=y$, has a singularity of the form $\sim d(x,y)^{-\mathrm{dim}_{\mathbf R} X+2}$ near $x=y$ ($G(x,y)$ has log pole if dim$_{\mathbf R} X = 2$), and is bounded from below. The infimum of the Green's function is itself an important invariant of the metric.
However, both the diameter and the Green's function depend on the metric in an intricate way. 

\bigskip
This dependence takes on a particular importance when the metric varies, and we consider the Gromov-Hausdorff convergence of sequences of metrics. A case in point is the following version of a theorem of Gromov:

\bigskip
{\bf \small Gromov compactness theorem}: Let $(X_\ell,g_\ell)$ be a sequence of compact Riemannian manifolds. If the sequence admits a uniform lower bound on their Ricci curvature and a uniform upper bound on their diameters, then it admits a subsequence which converges in the sense of Gromov-Hausdorff.

\bigskip
It is usually not easy to determine when a sequence of Riemannian metrics admits a uniform upper bound on their diameters. A classic result is the following

\bigskip
{\bf \small Bonnet-Myers theorem}: Let $(X,g_{ij})$ be a complete Riemannian manifold
of real dimension $m$. If it admits the following positive lower bound on its Ricci curvature
\bea
R_{ij}\geq (m-1)\,K\,g_{ij}
\eea
for some constant $K>0$, then it is compact, and the following bound for the diameter holds,
\bea
{\rm diam}_g\,X\leq {\pi\over\sqrt K}.
\eea

\bigskip
Similarly, we only have the following theorem of Cheng-Li   \cite{CL} on lower bounds for the Green's function on compact Riemannian manifolds:

\bigskip
{\bf\small Cheng-Li theorem}: Let $(X,g_{kj})$ be a compact Riemannian manifold of real dimension $m$. If its Ricci curvature satisfies the lower bound
\bea
R_{kj}\geq - (m-1)Kg_{kj}
\eea
for some positive constant $K$, then 
\bea
G(x,y)\geq -C{({\rm diam}_gX)^2\over V_g}
\eea
where $C$ is a constant depending only on $m$ and $K$.

\bigskip
We note the presence of a lower bound assumption on the Ricci curvature in all the above theorems. This presents a severe difficulty in many applications, notably on the Gromov-Hausdorff convergence of geometric flows, where it is usually hard to obtain bounds for the Ricci curvature.

\bigskip

\subsection{\bf THE K\"AHLER CASE}\label{section 2.3}
\hfill\break

This is where great progress on these issues has been made in recent years. We describe it now. As the crucial feature is the uniformity of the estimates with respect to the underlying metric, we have to specify carefully the classes of K\"ahler metrics we shall consider.

\bigskip
Thus let $(X,\o_X)$ be an $n$-dimensional compact K\"ahler manifold, equipped with a reference K\"ahler metric $\o_X$, normalized so that $V_{\o_X}=\int_X\o_X^n=1$.
Given a K\"ahler metric $\omega$ on $X$, recall that we denote by $V_\o= \int_X \omega^n$ its volume. We
define the relative volume function $F_\o$ by
$$
e^{F_\omega} = \frac{1}{V_\omega}\frac{\omega^n}{\omega_X^n}.
$$ 
Given $p\ge 1$ we define the $p$-th Nash-Yau entropy of $\omega$ by
$${\mathcal N}_p(\omega) = \frac{1}{V_\omega}\int_X \Big|\log \frac{1}{V_\omega} \frac{ \omega^n}{\omega_X^n} \Big|^p \;{\omega^n} =  \int_X |F_\omega|^p e^{F_\omega} \omega_X^n.$$
For given parameters $0< A\le +\infty, K>0$, we consider the following subsets of the space of K\"ahler metrics on $X$:
\begin{equation}\label{eqn:the set}
{\mathcal W}(n,p, A, K): = \Big\{ \omega :\; [\omega]\cdot [\omega_X]^{n-1}< A,\, {\mathcal N}_p(\omega) \le K \Big\}.
\end{equation}
The invariant $[\omega]\cdot[\omega_X]^{n-1}$ can be recognized as the intersection of the K\"ahler classes $[\omega]$ and $[\omega_X]$, and it is convenient to denote it by $I_\omega$,
\begin{equation}\label{norinter}
I_\omega = [\omega]\cdot[\omega_X]^{n-1}.
\end{equation}
Note that we allow subsets of K\"ahler classes with unbounded intersection numbers, and such admissible sets of K\"ahler metrics correspond to $A = \infty$ in \mref{eqn:the set} and are denoted by ${\mathcal W}(n, p, \infty, K)$. 

\bigskip
Then we have the following theorems, obtained by the authors in joint work with Jian Song and Jacob Sturm. The first one gives uniform estimates for the Green's function, the diameter, and the local volume growth for K\"ahler metrics in ${\mathcal W}(n,p,A,K)$:

\bigskip
\begin{theorem}\label{thm:gpss}
\cite{GPSS}
Let $\o$ be a K\"ahler metric in $\mathcal W(n,p,A,K)$. Then there are constants $C_0,C_1,C_2,C_3$ depending only on $X,\o_X,n,p,A,K$ and $\epsilon>0, \epsilon'>0$ depending only on $n$, and $p$ so that the following hold:

{\rm(a)}
\bea
{\rm inf}_{x,y}G_\o(x,y)\geq -{C_1\over V_\o};
\eea

{\rm (b)} We have the following integral bounds  
\bea
{\rm sup}_{x\in X}\big( V_\omega^\epsilon\int_X|G_\o(x,\cdot)|^{1+\epsilon}\o^n
+
V_\omega ^{\epsilon'}\int_X|\na G_\o(x,\cdot)|^{1+\epsilon'}\o^n\big)
\leq C_1
\eea

{\rm(d)} ${\rm diam}_\o X\leq C_2$;

{\rm(e)} There exists $\alpha = \alpha(n,p)>0$ such that, for any $R$ with $0<R<1$, we have
\bea
{{\rm Vol}_\o B_\o(x,R)\over V_\o}\geq {1\over C_3}R^\alpha
\eea
for any $x\in R$. Here $B_\o(x,R)$ denotes the ball centered at $x$ and of radius $R$ with respect to $\o$.
\end{theorem}

\bigskip
In the earliest version of this theorem, the classes of metrics considered had to satisfy an additional constraint, more specifically their volume forms had to be being bounded from below by a nonnegative function $\gamma$ whose vanishing locus is small in a suitable sense. This constraint was harmless in all applications considered. However, it turns out to be removable altogether, by Vu \cite{V} for the diameter estimates, and independently by Guedj and T\^o \cite{GT} and by the authors of \cite{GPSS} in \cite{GPSS24} for both the diameter and the Green's function. The argument in \cite{GPSS24} is just a simple modification of the original arguments in \cite{GPSS}. The above version is the later version without the constraint.

\bigskip

We stress that the classes ${\mathcal W}(n,p,A,K)$ of K\"ahler metrics are defined only by bounds on the Nash-Yau entropy ${\mathcal N}$ and the intersection number $I_\o=[\o]\cdot[\o_X]^{n-1}$. Thus the above bounds do not require any assumption on the Ricci curvature. If we keep in mind the fact that the Nash-Yau entropy is an integral of the relative volume form, and that the Ricci curvature requires two derivatives of this volume form (\ref{eqn:ric}), the above estimates have gained practically two derivatives. As an immediate consequence, we obtain the following powerful version for K\"ahler geometry of the Gromov compactness theorem:

\medskip
\begin{theorem}\cite{GPSS}
Let $(X,\o_X)$ be a fixed compact K\"ahler manifold, and let $(X,\o_\ell)$ be a sequence of K\"ahler manifolds in ${\mathcal W}(n,p,A,K)$. Then there exists a subsequence, still denoted by $(X,\o_\ell)$ for simplicity, such that
\bea
(X,\o_\ell)\to (Z,d_Z)
\eea
in the sense of Gromov-Hausdorff, where $(Z,d_Z)$ is a compact metric space.
\end{theorem}

\bigskip

The above progresses suffice to settle some long-standing questions about the limiting behavior of diameters in the K\"ahler-Ricci flow. Let $[0,T)$ be its maximum time interval of existence of the K\"ahler-Ricci flow \eqref{eqn:krf}.

\bigskip
\begin{theorem}\cite{GPSS}
Let $(X,\o_0)$ be a compact K\"ahler manifold.

{\rm (a)} Consider the K\"ahler-Ricci flow $t\to \o(t)$ defined by
\bea\label{eqn:krf}
\dot\o(t)=-Ric(\o),
\qquad
\o(0)=\o_0.
\eea
Assume that $T<\infty$ and ${\rm lim}_{t\to T^-}Vol(X,\o(t))> 0$. Then 
\bea
{\rm diam}_{\o(t)}X\leq C, \qquad t\in [0,T)
\eea
for some constant $C$ depending only on $n$ and $\o_0$.

\smallskip
{\rm (b)} Consider the normalized K\"ahler-Ricci flow defined by
\bea
\dot\o(t)=-Ric(\o(t))-\o(t),
\qquad \o(0)=\o_0.
\eea
Assume that $T=\infty$, then 
\bea
{\rm diam}_{\o(t)}X\leq C,\qquad t\in [0,\infty).
\eea
\end{theorem}

\bigskip

As shown in further joint work of the authors with Jian Song and Jacob Sturm,
the analysis can be pushed much farther, giving Sobolev-type inequalities as well as other basic inequalities for geometric analysis such as estimates for the heat kernel and for the eigenvalues of the Laplacian. Thus we also obtain the following Sobolev-type inequality:

\bigskip

\begin{theorem}\label{thm:Sob} \cite{GPSS23, GPSS24}
Given $p>n$ and $K>0$, there exist   $q = q(n, p)>1$ and $C = C(n, p,   K, q)>0$ such that  for any K\"ahler metric $\omega\in {\mathcal W}(n,p,\infty, K )$ and any $u\in W^{1,2}(X)$,  we have the following Sobolev-type inequality 
\begin{equation}\label{eqn:Sob}
\Big(\frac{1}{V_\omega}\int_X | u - \overline{u}  |^{2q}\omega^n   \Big)^{1/q}\le C \frac{ I_\omega}{V_\omega} \int_X |\nabla u|_\omega^2 \omega^n,
\end{equation}
where $\overline{u} = \frac{1}{V_\omega}\int_X u \omega^n$ is the average of $u$ over $(X,\omega)$.
\end{theorem}

\bigskip
We remark that this Sobolev inequality \mref{eqn:Sob} is {\em  scale invariant}, although the exponent $q$ may be  smaller than $\frac{n}{n-1}$. If we impose a stronger condition on $e^{F_\omega}$, i.e. its $L^{1+\epsilon'}(X,\omega_X^n)$-norm is uniformly bounded for some $\epsilon'>0$, then the constant $q>1$ can be chosen as close as possible  to the {exponent} $\frac{n}{n-1}$ as in the Euclidean case, with the constant $C_q$ possibly blowing up as $q$ approaches $ \frac{n}{n-1}$. However, for most applications, for example Moser iteration, the exponent $q>1$ suffices. It is also remarkable that the analogue of the usual Sobolev constant now depends only on $n,p,K,q, V_\o$ and $I_\o$.
For lack of space, we do not include here the estimates for the heat kernel and the eigenvalues of the Laplacian. They can be found in \cite{GPSS23}. A stronger version, valid for K\"ahler spaces with suitable singularities, is described below.

\bigskip

\subsection{GEOMETRIC INEQUALITIES ON K\"AHLER SPACES}
\hfill\break

It is a very important property of the geometric inequalities which we have discussed so far that they descend on suitable K\"ahler spaces. This is a property which is even stronger than their uniformity with respect to degenerating families of smooth K\"ahler manifolds, as it opens the door to analysis on singular K\"ahler spaces. 

\bigskip
We need to specify the type of singular K\"ahler spaces which will be considered. A key requirement is that the corresponding classes of K\"ahler metrics can include singular K\"ahler-Einstein metrics. With this in mind, we introduce the following K\"ahler spaces and their corresponding classes of K\"ahler metrics.

\bigskip
Let $X$ be an $n$-dimensional compact normal K\"ahler space, and let $Y$, $\pi: Y\to X$  be a nonsingular model of $X$. It is known that $Y$ is a smooth K\"ahler manifold (c.f. Lemma 2.2 \cite{CMM}). Fix then a smooth K\"ahler metric $\theta_Y$ on $Y$. We can then define the space of semi-K\"ahler currents
$$
{\mathcal {AK}}(X,\theta_Y,n,p,A,K)
$$
by the following conditions:

\bigskip

\begin{enumerate}

\item $[\omega]$ is a K\"ahler class on $X$ and $\omega$ has bounded local potentials. 

\smallskip

\item $[\pi^*\omega] \cdot [\theta_Y]^{n-1}\leq A$. 

\smallskip

\item The $p$-th Nash-Yau entropy is bounded for some $p>n$, i.e.
$${\mathcal{N}}_p (\omega) = \frac{1}{V_\omega} \int_Y \left|\log \frac{1}{V_\omega} \frac{(\pi^*\omega)^n}{(\theta_Y)^n} \right|^p (\pi^*\omega)^n  \leq K, $$
where $V_\omega= [\omega]^n$.

\medskip

\item The log volume measure ratio 
$$\log \left( \frac{(\pi^*\omega)^n}{(\theta_Y)^n} \right)$$ 
has log-type analytic singularities. 

%
\end{enumerate}

\bigskip

Recall that a function $f$ on $Y$ is said to have log-type analytic singularities if there exist holomorphic effective divisors $D_j$, $1\leq j\leq N$, with simple normal crossings such that
\bea
f=\sum_{k=1}^Ka_k(-\log)^k(\p_{j=1}^N e^{f_{k,j}} |\sigma|_{h_j}^{2b_{k,j}})
\eea
where $a_k\in{\bf R}$, $b_{k,j}\geq 0$, $f_{k,j}\in C^\infty(Y)$, $K$ is a positive integer, and $\sigma_k$is a defining section for $D_j$, and $h_j$ is a smooth Hermitian metric associated to $D_j$.
We note that the condition that the volume measure ratio has log type analytic singularities is meant to guarantee the smoothness of $\omega$ in an open Zariski subset of $X$. 

\bigskip
For any $\omega \in \mathcal{AK}(X, \theta_Y, n, p, A, K)$, we let $\cS_{X, \omega}$ be the union of the singular set of $X$ and the singular set of $\pi_*\left( \log \frac{(\pi^*\omega)^n}{(\theta_Y)^n} \right)$. By definition, $\cS_{X, \omega}$ is an analytic subvariety of $X$. The current $\omega$ is a smooth K\"ahler metric on $X\setminus \cS_{X, \omega}$ (c.f. Proposition 7.1 in \cite{GPSS23}). On $X\setminus \cS_{X,\o}$, we can introduce the distance function $d$
\bea
d(x,y)={\rm inf}\{L_\o(c)=\int_0^1|\dot c(t)|_\o dt; \ c:[0,1]\to X\setminus\cS_{X,\o},
\ c(0)=x,\quad c(1)=y\}\nonumber
\eea
where $c$ ranges over piecewise differentiable curves. The space $(X\setminus\cS_{X,\o},d)$ is then a metric space. Let $(\hat X,d)$ be its metric completion,
\bea
(\hat X, d) = \overline{(X\setminus \cS_{X, \omega}, \omega|_{X\setminus \cS_{X, \omega}})}.
\eea
Finally we define the Sobolev space $W^{1,2}(\hat X,d,\o^n)$ as the space of all $u: X\to \bf R$ such that $u_{\vert_{\mathcal K}}\in W^{1,2}({\mathcal K})$ for all ${\mathcal K}\subset\subset X\setminus\cS_{X,\o}$, and
\bea
{\rm sup}_{{\mathcal K}\subset\subset X\setminus\cS_{X,\o}}\|u\|_{W^{1,2}({\mathcal K})}
=
{\rm sup}_{{\mathcal K}\subset\subset X\setminus\cS_{X,\o}}
(\int_{\mathcal K}(|u|^2+|\na u|_\o^2)\o^n)^{1\over 2}<\infty.
\eea
We can then define the $W^{1,2}$ norm of $u$ to be
\bea
\|u\|_{W^{1,2}}^2={\rm sup}_{{\mathcal K}\subset\subset X\setminus\cS_{X,\o}}
\|u\|_{W^{1,2}({\mathcal K})}^2.
\eea
We have then

\begin{theorem} \label{singsob} Let $X$ be an $n$-dimensional compact normal K\"ahler space. For any $\omega\in \mathcal{AK}(X, \theta_Y, n, p, A, K)$, the metric measure space $(\hat X, d, \omega^n)$ associated to $(X, \omega)$ satisfies the following properties. 

\begin{enumerate}

\item There exists $C=C(X, \theta_Y, n, p, A, K)>0$ such that  
$${\textnormal{diam}}(\hat X, d) \leq C.$$
In particular, $(\hat X, d)$ is a compact metric space.

\medskip 

\item  There exist $q>1$ and $C_S=C_S(X, \theta_Y, n, p,  A, K,  q)>0$ such that 
$$
\Big(\frac{1}{V_\omega}\int_{\hat X} | u  |^{2q}\omega^n   \Big)^{1/q}\le \frac{C_S}{V_\omega} \left( \int_{\hat X} |\nabla u|^2 ~\omega^n + \int_{\hat X} u^2 \omega^n \right) .
$$
for all $u\in W^{1, 2}(\hat X, d, \omega^n)$. 

\medskip

\item There exists $C=C(X, \theta_Y, n, p, A, K)>0$ such that the following trace formula holds for the heat kernel of $(\hat X, d, \omega^n)$
$$H(x,x, t) \leq \frac{1}{V_\omega} + \frac{C}{V_\omega} t^{-\frac{q}{q-1}}. $$ 

\medskip

\item Let $0=\lambda_0 < \lambda_1 \leq \lambda_2 \leq ... $ be the increasing sequence of eigenvalues of the Laplacian $-\Delta_\omega$ on $(\hat X, d, \omega^n)$. Then there exists $c=c(X, \theta_Y, n, p, A, K)>0$ such that
$$\lambda_k \geq c k^{\frac{q-1}{q}}. $$

\end{enumerate}

\end{theorem}

Besides the ideas which we have described, an important ingredient of the proof in the case of K\"ahler spaces with singularities is the cut-off functions constructed by Sturm \cite{Stu}.

\medskip

As stressed in \cite{GPSS23},  
Theorem \ref{singsob} appears to be the first to provide general formulas for the Sobolev inequality, the heat kernel,  and eigenvalue estimates on complex spaces with singularities. Such estimates were only known  previously to hold in special cases such as the following:
\begin{enumerate}
\item $X$ is a normal complex space and $\omega$ is a Bergman metric via a global projective embedding or the restriction of a smooth metric via local Euclidean embeddings (cf. \cite{LT, Be}).

\smallskip

\item $(X, \omega)$ is a sequential Gromov-Hausdorff limit of K\"ahler-Einstein manifolds (cf. \cite{DS}). 

\end{enumerate}

\bigskip

\subsection{\bf THE PROOF: THREE BASIC IDEAS} 
\hfill\break

The proof of all the preceding results rests on the following three basic ideas (a),(b), and (c).

\bigskip
{\bf (a) The K\"ahler condition and the Monge-Amp\`ere equation}

\medskip
We have seen earlier how the Calabi conjecture can be reduced to a complex Monge-Amp\`ere equation. The basic underlying idea was that a K\"ahler metric can be determined by a complex Monge-Amp\`ere equation from its cohomology class and its volume form. However, there is freedom in the choice of reference metric within a given K\"ahler class, and in reducing this volume form to a scalar function. Thus let $X$ be a compact $n$-dimensional K\"ahler manifold, and fix a K\"ahler reference metric $\o_X$ normalized so that $V_X=\int_X\o_X^n=1$. If $\o$ is any other K\"ahler metric, and we choose a closed smooth $(1,1)$-form $\theta\in [\o]$, then we can write $\o=\theta+i\ddb\varphi$ by the $\ddb$-Lemma. If we define the scalar volume function of $\o$ by
\bea
F_\o=\log ({1\over V_\o}\cdot{\o^n\over\o_X^n}),
\eea
then it is just a tautology that $\varphi$ satisfies
\bea
(\theta+i\ddb\varphi)^n=V_\o e^{F_\o}\o_X^n,\qquad
\theta+i\ddb\varphi>0.
\eea

Nevertheless, we can view this tautology as a complex Monge-Amp\`ere equation for $\varphi$.
This simple shift from using $\o_X^n$ instead of $\theta^n$ as a reference volume form  is natural if we want to let $[\o]$ degenerate. More important, the Nash-Yau entropy is just the $L(\log L)^p$ norm of the right hand side. Also, if $\o$ is assumed to be in one of the classes ${\mathcal W}(X,\o_X,n,A,K)$, and we would like uniform estimates with respect to this class, it is useful to choose some specific representative $\theta\in [\o]$. In \cite{GPSS23}, it was shown that the choice of $\theta$ which is harmonic with respect to the metric $\o_X$ gives a form $\theta$ which is uniformly bounded in $C^3$ norm.
So the main issue is to determine how much of the geometry of $\o$ can be recaptured in this set-up from the Nash-Yau entropy. A broad goal of getting geometric estimates, rather than $C^k$ estimates, may have appeared first in print in a 2020 paper by X. Fu, B. Guo, and J. Song \cite{FGS}. However, at that time, there were as yet no tool sufficiently powerful for the kind of estimates that we are aiming for.

\bigskip
{\bf (b) The method of auxiliary Monge-Amp\`ere equations}

\medskip

If we turn for guidance to Yau's solution of the Calabi conjecture, we see that the $C^2$ estimates for $\varphi$, and hence the $C^0$ estimate for the corresponding metric $\o$, cannot be obtained from his estimates from just the $L^1(\log L)^p$ norm of the right hand side. Even the $L^\infty$ bound for $\varphi$ cannot be obtained from Yau's method of Moser iteration, since it requires assuming an $L^q$ bound on the right hand side for some $q>n$.

\medskip
But at least for the $L^\infty$ estimates for $\varphi$, we have the sharp result of Kolodziej \cite{K}, obtained some 20 years after Yau's original solution, that an $L^\infty$ estimate for $\varphi$ can indeed be obtained assuming a bound on the $L^1(\log L)^p$ norm of the right hand side for some $p>n$. Kolodziej's result is very encouraging. Nevertheless, geometric estimates are defined in terms of the metric $\o$, and naively, it seems that $C^2$ estimates for $\varphi$ cannot be bypassed.

\medskip
Thus we do seem to need a new method. A good test for such a method is that it should at least be able to reproduce Kolodziej's result, which had been established using pluripotential theory. This had been in itself an open problem since the appearance of Kolodziej's result, and it was solved in 2021 by the authors, in joint work with Freid Tong \cite{GPT}. The basic idea is a comparison to an auxiliary complex Monge-Amp\`ere equation, whose solution exists by Yau's theorem. Now comparisons with another PDE is a well-known method in PDE theory, and valuable comparisons with a Monge-Amp\`ere equation had been used by Dinew and Kolodziej \cite{DK} and more recently by Chen and Cheng \cite{CC} in their study of the constant scalar curvature equation. The key in the present situation is the choice of the auxiliary Monge-Amp\`ere equation and of the comparison function. 

\medskip

\medskip
$\bullet$ Let us assume that the solution $\varphi$ of the equation that we are considering has been normalized to ${\rm sup}_X\varphi=0$. We would like to show that
$\varphi$ is bounded from below, or equivalently, that the set 
$$
\Omega_s:=\{-\varphi-s>0\}
$$
is empty for some constant $s=S_\infty$. For this we use the classical method of De Giorgi, which was also the method used by Kolodziej. In this method, we have to find a continuous, non-negative, and decreasing function $\phi:{\bf R}\to {\bf R}_{\geq 0}$ with the property that

\begin{enumerate}

\item $\varphi\geq -S_\infty$ if and only if $\phi(s)=0$ for $s\geq S_\infty$;

\item The existence of such an $S_\infty$ would be guaranteed by the existence of $\delta>0$ so that an inequality of the following form holds
\bea
\label{DeG}
r\phi(s+r)\leq C_\delta \phi(s)^{1+\delta},
\qquad r,s>0.
\eea
\end{enumerate}

\bigskip
In De Giorgi's original approach to the De Giorgi-Nash-Moser theory, the function $\phi(s)$ was chosen to be
$$
\phi(s)={\rm Volume}(\Omega_s).
$$
In Kolodziej's proof of sharp $L^\infty$ estimates for the complex Monge-Amp\`ere equation, $\phi(s)$ was chosen to be
$$
\phi(s)={\rm Capacity}(\Omega_s).
$$
In our new approach to $L^\infty$ estimates, we shall choose
$$
\phi(s)=\int_{\Omega_s}e^{F_\o}\omega_X^n.
$$

\medskip

The main task now is to show that $\phi(s)$ satisfies a growth condition (\ref{DeG}) which insures that it vanishes beyond a certain $s=S_\infty$. To do this, we introduce the
function $A_s$ 
\bea
\label{As}
A_s=  \int_{\Omega_s}(-\varphi-s)e^{F_\o}\o_X^n.
\eea
Since trivially $\Omega_{s+r}\subset \Omega_s$ for any $r>0$, we have
\bea
A_s\geq   \int_{\Omega_{s+r}}(-\varphi-s)e^{F_\o}\o_X^n
\geq
r   \int_{\Omega_{s+r}}e^{F_\o}\o_X^n
=r\phi(s+r),
\eea
the desired growth condition for $\phi(s)$ would immediately follow if we can show that
\bea
A_s\leq C_\delta\,\phi(s)^{1+\delta}
\eea
for some constant $C_\delta>0$. We observe that this is reminiscent of a reverse H\"older inequality, which can usually be established only for solutions of an elliptic partial differential equation.

\bigskip
$\bullet$ The key now is the choice of an appropriate auxiliary complex Monge-Amp\`ere equation. To illustrate the main idea and simplify the notations, we may assume $\theta$ is a {\em nonnegative} $(1,1)$-form, otherwise the $L^\infty$ estimate should be modified by the envelope of $\theta$.  One natural candidate, in view of the previous discussion, is
\bea
(\theta+i\ddb\psi_s)^n
=V_\o{(-\varphi-s)_+\over A_s}e^{F_\o}\o_X^n,
\qquad \theta+i\ddb\psi_s>0,
\eea
where $(-\varphi-s)_+$ denotes the positive part of the function $-\varphi-s$.  Note that the constant $A_s$, as defined previously in (\ref{As}), is precisely the normalization constant for both sides in the above equation to have the same integral over $X$.
However, 
since the function $(-\varphi-s)_+$ is not smooth, we approximate it by smooth functions $\tau_k(-\varphi-s)$, where $0<\tau_k(t)\downarrow t_+$ is a sequence of smooth positive functions, and consider the auxiliary Monge-Amp\`ere equation
\bea
(\theta+i\ddb\psi_{s,k})^n=
V_\o{\tau_k(-\varphi-s)\over A_{s,k}}e^{F_\o}\o_X^n,
\qquad \theta+i\ddb\psi_{s,k}>0,
\eea
where the normalization constant $A_{s,k}$ is now defined by
\bea
A_{s,k}=  \int_X\tau_k(-\varphi-s)e^{F_\o}\o_X^n.
\eea
Yau's theorem applies now, and guarantees the existence of a unique smooth solution $\psi_{s,k}$ to the auxiliary equation satisfying the normalization ${\rm sup}_X\psi_{s,k}=0$.

\bigskip
$\bullet$ So far, our discussion of the strategy for the auxiliary Monge-Amp\`ere equation has been completely general, and we have not even specified the equation for $\varphi$ that we are considering. Assume now that $\varphi$ satisfies the complex Monge-Amp\`ere equation (\ref{MA1}). Then we can show using the maximum principle, that we have the following comparison inequality for $-\varphi-s$,
\bea
\label{comparison}
{-\varphi-s\over A_{s,k}^{1\over n+1}}\leq \epsilon_n(-\psi_{s,k}+A_{s,k})^{n\over n+1}
\eea
where $\epsilon_n$ is a suitable constant depending only on $n$.

\bigskip
$\bullet$ One concern about the choice of auxiliary Monge-Amp\`ere equation may have been that it and its solution $\psi_{s,k}$ depend themselves on the unknown function $\varphi$. But it turns out that we need to know very little about it, besides the fact that $\psi_{s,k}$ is plurisubharmonic with respect to $\theta$. Thus we note that the comparison inequality (\ref{comparison}) implies
\bea
\int_{\Omega_s}{\rm exp}
\big\{\beta_0{(-\varphi-s)^{n+1\over n}\over A_{s,k}^{1\over n}}\big\}\o_X^n
\leq
{\rm exp}(C_n\beta_0A_{s,k})\int_{\Omega_s}
{\rm exp}(-C_n\beta_0\psi_{s,k})\o_X^n,
\eea
with $C_n=\epsilon_n^{n+1\over n}$. By a classic inequality of H\"ormander \cite{Ho} for plurisubharmonic functions, extended to the global setting by Tian \cite{Ti}, for a suitable constant $\beta_0$ small enough and all $\theta\leq \kappa\,\o_X$, the integral on the right hand side is uniformly bounded by a constant depending only on $n$, $\beta_0$, and $\kappa$. We can now take the limit as $k\to\infty$ and obtain
\bea
\label{exp}
\int_{\Omega_s}
{\rm exp}\big\{\beta_0{(-\varphi-s)^{n+1\over n}\over A_s^{1\over n}}\big\}\o_X^n
\leq
K_{n,\beta_0,\kappa}{\rm exp}(C_n\beta_0A_s)
\eea
for a constant $K_{n,\beta,\kappa}$ depending only on $n,\beta_0,\kappa$.

\bigskip
$\bullet$ It is now easy to derive an estimate for $A_s$ from this inequality, since $A_s$ is an integral of $-\varphi-s$ with respect to the measure $e^{F_\o}\o_X^n$, while the above left hand side is the integral of its exponential with respect to the measure $\o_X^n$. For this, we apply Young's inequality, which leads to, after some simplifications and using the key inequality (\ref{exp}),
\bea
\int_{\Omega_s}(-\varphi-s)^{p(n+1)\over n}e^{F_\o}\o_X^n
\leq C\,A_s^{p\over n}
\eea
where $C$ is a constant depending only on $X,\o_X,n,p,\|e^{F_\o}\|_{L^1(\log L)^p}$ for any $p\geq 1$. Note the appearance here of the Nash-Yau entropy $\|e^{F_\o}\|_{L^1(\log L)^p}$. From here, we obtain by the H\"older inequality
\bea
A_s&=&\int_{\Omega_s}(-\varphi-s)e^{F_\o}\o_X^n
\leq(\int_{\Omega_s}(-\varphi-s)^{p(n+1)\over n}e^{F_\o}\o_X^n)^{n\over (n+1)p}
(\int_{\Omega_s}e^{F_\o}\o_X^n)^{1\over q}
\nonumber\\
&\leq& C\,A_s^{1\over n+1}\phi(s)^{1\over q}
\nonumber
\eea
where $q$ is defined by ${n\over p(n+1)}+{1\over q}=1$, and hence
\bea
A_s\leq C\,\phi(s)^{1+n\over qn}.
\eea
The exponent ${1+n\over qn}$ works out to be $1+\delta_0$, with $\delta_0={p-n\over pn}$. Thus we have $\delta_0>0$ and the desired growth rate for $\phi(s)$ if the Nash-Yau entropy $\|e^{F_\o}\|_{L^1(\log L)^p}$ is bounded for some $p>n$. This is the sharp $L^\infty$ estimate obtained previously by Kolodziej using pluripotential theory, in a form which encompasses the subsequent uniform versions obtained by Eyssidieux, Guedj, Zeriahi \cite{EGZ} and Demailly and Pali \cite{DP}.

\bigskip
It was already shown in \cite{GPT} that the proof there for Monge-Amp\`ere equations applied equally well to fully non-linear equations satisfying a simple structural condition. This condition was shown by Harvey and Lawson \cite{HL,HL1} to be satisfied for large classes of equations, including all invariant G{$\mathring{\text{a}}$}rding-Dirichlet equations. In the case of fully non-linear equations, it is useful to have a uniform bound of the energy by the entropy. Such a bound is supplied in \cite{GPb}, also by the method of auxiliary Monge-Amp\`ere equations. The above $L^\infty$ estimates can be even extended to fully non-linear equations in Hermitian geometry \cite{GPa}, if the auxiliary equation is replaced by the Dirichlet problem for the complex Monge-Amp\`ere equation, and the maximum principle is replaced by the ABP inequality, as pioneered by Blocki \cite{Bl} in his proof of $L^\infty$ estimates for the complex Monge-Amp\`ere equation assuming the more restrictive condition of finiteness of $\|e^{F_\o}\|_{L^{2+\epsilon}}$ for $\epsilon>0$. A brief survey of these developments has been given in \cite{GP}. Below we concentrate rather on geometric estimates for the Monge-Amp\`ere equation.

\bigskip
{\bf (c) Estimates for the Green's function}

\medskip
The third idea in this new approach to geometric estimates is to start from estimates for the Green's function and derive geometric estimates from them, instead of the other way around. Indeed, assume that we have integral inequalities for the Green's function of the form stated in Theorem \ref{thm:gpss}. Let $z_0,w_0$ be points on $X$ with $d_\o(z_0,w_0)={\rm diam}_{\o}X$. Set $\rho(z)=d_\o(z_0,z)$, which is a Lipschitz function satisfying $|\na\rho|\le 1$. Then Green's formula implies
\bea
0=\rho(z_0)={1\over V_\o}\int_X \rho\o^n
+
\int_X\<\na G(z_0,\cdot),\na \rho(\cdot)\>\o^n
\eea
and
\bea
{\rm diam}_\o X=\rho(w_0)={1\over V_\o}\int_X\rho\o^n
+
\int_X\<\na G(w_0),\cdot),\na\rho(\cdot)\>\o^n.
\eea
Subtracting the two formulas gives
\bea
{\rm diam}_\o X&=&\int_X\<\na G(w_0,\cdot)-\na G(z_0,\cdot),\na\rho(\cdot)\>\o^n\nonumber\\
&\leq&
2\,{\rm sup}_{z\in X}\int_X |\na G(z,\cdot)|\o^n
=
2\,{\rm sup}_{z\in X}\|\na G(z,\cdot)\|_{L^1}.
\eea
Similarly, if $0<R<1$ and we introduce a cut-off function $\eta\in C_0^1(B(z_0,R))$ with $\eta=1$ in $B(z_0,R/2)$  and $|\na\rho|< 4R^{-1}$, Green's formula can be applied to the function $\eta\,\rho$. The same argument together with H\"older's inequality would give then a lower bound of the form $R^\alpha$ for ${\rm Vol}_\o B(z_0,R)$ and some strictly positive exponent $\alpha$, if we can assume a bound on 
${\rm sup}_{z\in X}\|\na G(z,\cdot)\|_{L^{1+\epsilon}}$ for some fixed $\epsilon>0$.

\bigskip
The switch to estimates for the Green's function revealed its advantages when the above method of auxiliary Monge-Amp\`ere equations turns out to apply equally well to the Laplace equation, and not just fully non-linear equations. A model estimate which can be obtained in this manner is the following mean-value inequality:

\bigskip
\begin{theorem} \cite{GPS, GPSS}
Let $(X,\o_X)$ be a compact K\"ahler manifold. Let $\o\in {\mathcal W}(n,p,A,K)$ be a K\"ahler metric, and let $v\in L^1(X)$ satisfy $\int_X v\o^n=0$ and 
\bea
v\in C^2(\bar\Omega_0),\quad \Delta_\o v\geq -1\ \ {\rm in}\ \ \Omega_0 = \{v>0\}.
\eea
Then we have
\bea \label{eqn:2.60}
{\rm sup}_Xv\leq C\,(1+{1\over V_\o}\int_X |v|\o^n)
\eea
for some constant $C$ depending only on $\o_X,n,p,K$.
\end{theorem}

\bigskip
\noindent
{\it Proof}. 
We only mention some important ingredients of the proof.
We may assume that $\|v\|_{L^1(X,\o^n)}\leq 1$, otherwise we replace $v$ by $vV_\o/\|v\|_{L^1(X,\o^n)}$. It suffices then to show that ${\rm sup}_Xv$ is bounded by a constant.
We again rely on an auxiliary complex Monge-Amp\`ere equation, in this case, the following equation
\bea
(\theta+i\ddb\psi_{s,k})^n=V_\o{\tau_k(v-s)\over A_{s,k}}e^{F_\o}\o_X^n,
\quad
\theta+i\ddb\psi_{s,k}>0,
\eea
where $A_{s,k}$ is defined by
\bea
A_{s,k}=\int_X\tau_k(v-s)e^{F_\o}\o_X^n.
\eea
If we define the comparison function $\Phi$ by
\bea
\Psi=-\epsilon(-\psi_{s,k}+\varphi+\Lambda)^{\frac{n}{n+1}}+v-s
\eea
we find then that $\Psi\leq 0$ for suitable constants $\epsilon$ and $\Lambda$. If we recall that the set $\Omega_s$ was defined by $\Omega_s=\{v>s\}$ and apply the De Giorgi iteration argument to the function $\phi(s)$ defined by
\bea
\phi(s)=\int_{\Omega_s}e^{F_\o}\o_X^n,
\eea
we can show as before that there exists a constant $S_\infty$ with $\phi(s)=0$ for $s>S_\infty$. Thus $v\leq S_\infty$, and the theorem is proved.

\bigskip
The preceding theorem was formulated so as to apply most readily to $v(y)=-V_\o\,G(x,y)$, for a fixed $x\in X$. We obtain in this manner
\bea
-{\rm inf}_{y\in X}G(x,y)\leq \frac{C}{V_\omega}(1+\|G(x,\cdot)\|_{L^1}).
\eea

\bigskip
To go further, we also need a version of the preceding theorem which implies that $|v|\leq C$ for all functions $v\in C^2(X)$ satisfying $|\Delta_\o v|\leq 1$ and $\int_X v\o^n=0$. This is obtained using again an auxiliary Monge-Amp\`ere equation, this time
\bea
(\o+i\ddb\psi)^n
=
{{\rm max}(v,0)\over B}\o^n,
\qquad \o+i\ddb\psi>0,
\eea
where $B$ is defined by the normalization condition
\bea
B={1\over V_\o}\int_{X} \max (v, 0)\o^n.
\eea
The function ${\rm max}(v,0)$ is not smooth, so strictly speaking, we need to approximate it by smooth positive functions, as we had done in previous applications of the auxiliary Monge-Amp\`ere equation method. We omit these technicalities. The comparison function that we need in the present case turns out to be
\bea
\Phi=-\epsilon(-\psi+\epsilon^{n+1})^{n\over n+1}+{\rm max}(v,0).
\eea
For a suitable constant $\epsilon$, we have $\Phi\leq 0$. The rest of the proof follows as before.

\bigskip
We do not have space here to discuss the proofs of the other theorems stated in Section \S \ref{section 2.3}. Nevertheless, we would like to mention one more idea, which has not appeared as yet so far. It occurs in the proof of e.g, bounds for the integral
\bea
\int_X|G(x,\cdot)|^{1+\epsilon}\o^n.
\eea
Since $G(x,y)$ is bounded from below, we can view this integral as heuristically
\bea
\int_X {\mathcal G}(x,\cdot)^{1+\epsilon}\o^n
\eea
where ${\mathcal G}(x,y):=G(x,y)+C\geq 1$ for a suitable constant $C$. Let $H_k(y)=
{\rm min}({\mathcal G}(x,y),k)$, suitably smoothed out. We view the preceding integral as an approximation of the integral
\bea
\int_X{\mathcal G}(x,\cdot)H_k(\cdot)^\epsilon \o^n
\eea
and as such, it should be closely related to the solution $u_k$ of the equation
\bea
\Delta_\o u_k=-H_k^\epsilon+{1\over V_\o}\int_X H_k^\epsilon\,\o^n,
\quad
{1\over V_\o}\int_Xu_k\o^n=0.
\eea
The method of auxiliary Monge-Amp\`ere equations can now proceed in analogy with the proof of the mean-value inequality.

\bigskip

\section{\bf MORE RECENT DIRECTIONS PIONEERED BY YAU}
\setcounter{equation}{0}

As mentioned in the Introduction, this section is devoted to a survey of some of the directions pioneered by Yau since his 1976 landmark work on the Calabi conjecture. 
There are indeed many of them, with an abundance of deep and open problems which should be a fertile area for research in complex geometry and partial differential equations in the foreseeable feature. Our survey is necessarily brief, and we shall provide references to fuller surveys on each direction whenever possible.

\subsection{\bf THE DONALDSON-UHLENBECK-YAU THEOREM}

\hfil\break

Another celebrated work of Yau, some 10 years after his solution of the Calabi conjecture, is his solution with Karen Uhlenbeck of the Hermitian-Einstein equation.

\medskip
Let $(X,\o)$ be a compact K\"ahler manifold, and let $E\to X$ be a holomorphic vector bundle over $X$. A Hermitian-Einstein metric $H$ on $E$ is a Hermitian metric whose curvature form satisfies the following equation
\bea
\label{HE}
F\wedge {\o^{n-1}\over (n-1)!}=\lambda {\o^n\over n!}
\eea
for a constant $\lambda$, where $F\in \Lambda^{1,1}\otimes End(E)$ is the curvature form of $H$. In a local trivialization where the holomorphic coordinates on $X$ are given by $z^j$, $1\leq j\leq n$, and the sections of $E$ are given by $\varphi^\alpha$, $1\leq\alpha\leq r={\mathrm {rank}}\,E$, we can express the metric $H$ as $H=H_{\bar\alpha\beta}$, and the curvature form as 
\bea
F=F_{\bar kj}{}^\alpha{}_\beta
=
-\pl_{\bar k}(H^{\alpha\bar\gamma}\pl_jH_{\bar\gamma\beta}).
\eea
The Hermitian-Einstein equation (\ref{HE}) can then be written explicitly as
\bea
g^{j\bar k}F_{\bar kj}{}^\alpha{}_\beta=\lambda\,\delta^\alpha{}_\beta.
\eea
It can be viewed as a version of the Yang-Mills equation in the K\"ahler setting.
K\"ahler-Einstein metrics can be viewed as the particular case of Hermitian-Einstein metrics, when $E=T^{1,0}(X)$, and the Hermitian-Einstein metric coincides with the metric $\o$ on the base. Thus the K\"ahler-Einstein equation is more non-linear, but the Hermitian-Einstein equation requires solving for a whole metric, instead of just a scalar.
The following theorem was proved independently by Donaldson and Uhlenbeck-Yau:

\bigskip
{\bf\small Theorem of Donaldson-Uhlenbeck-Yau} \cite{Do, UhYa} Let $E\to X$ be a holomorphic vector bundle over a compact $n$-dimensional K\"ahler manifold $(X,\o)$. Assume for simplicity that $E$ is irreducible. Then $E\to X$ admits a Hermitian-Einstein metric if and only if $E\to X$ is stable in the sense of Mumford, and $\lambda$ is given by 
\bea
\lambda={1\over [\o^n]}\int_Xc_1(F)\wedge \o^n.
\eea
where $F_S$ is the curvature of any Hermitian metric on $S$.

\medskip
Recall that the bundle $E\to X$ is said to be stable in the sense of Mumford if, for any subsheaf ${\mathcal S}$ of $E$, we have
\bea
\mu({\mathcal S}) < \mu(E)
\eea
unless ${\mathcal S}={\mathcal O}(E)$. Here $\mu({\mathcal S})$ is the slope of ${\mathcal S}$, defined by
\bea
\mu({\mathcal S})= \frac{\mathrm{deg}_\omega (\mathcal S)}{\mathrm{rank}(\mathcal S)}.
\eea

\bigskip

The approaches of Donaldson and of Uhlenbeck-Yau are completely different. The one of Uhlenbeck-Yau may be viewed as a ``pure PDE" proof, as it does not rely on the highly non-trivial theorem of Mehta and Ramanathan \cite{MeRa1, MeRa2} that a stable bundle restricted to a generic hypersurface is stable, as the proof of Donaldson does. Rather, assuming that the bundle $E$ does not admit a Hermitian-Einstein metric, it shows by pure PDE methods that a destabilizing sheaf must exist. As the Uhlenbeck-Yau method is very powerful, and appears not to have been used since in other contexts, it may be appropriate to sketch it here.

\medskip
Fix a reference Hermitian metric $(H_0)_{\bar\alpha\beta}$ on $E$. Then any other metric $H_{\bar\alpha\beta}$ can be identified with the corresponding endomorphism $h$,
\bea
h^\gamma{}_\beta=(H_0)^{\gamma\bar\alpha}H_{\bar\alpha\beta}.
\eea
which is a positive endomorphism with respect to both $H_0$ and $H$ of $E$. It is convenient to introduce the notation
\bea
(\Lambda F)^\alpha{}_\beta=g^{j\bar k}F_{\bar kj}{}^\alpha{}_\beta.
\eea
Thus the Hermitian-Einstein equation can be written as
\bea
\Lambda F-\lambda I=0.
\eea

\bigskip
In the Uhlenbeck-Yau proof, one introduces the approximate Hermitian-Einstein equation
\bea
\Lambda F-\lambda I=-\epsilon \log h
\eea
and shows that it always admits a solution for $\epsilon>0$. Furthermore, by a suitable choice of reference metric $H_0$, the corresponding solution $h_\epsilon$ always satisfies ${\rm det}\,h_\epsilon =1$. The issue is then to determine when the sequence $\{h_\epsilon\}$ admits a subsequence converging in $C^\infty$, in which case the limit is then the desired Hermitian-Einstein metric. 

\medskip
$\bullet$ The first step is a reduction to a $C^0$ estimate. More precisely, if the endomorphisms $h_\epsilon$ are uniformly bounded, which is equivalent in view of the fact that they are positive to a uniform bound for their traces,
\bea
{\rm Tr}\,h_\epsilon \leq C<\infty
\eea
then a priori estimates to all orders can be derived that show that such a subsequence must exist.

\medskip
$\bullet$ One can proceed now by contradiction. Assume that the above inequality is violated for a subsequence, still denoted $\{h_\epsilon\}$ for simplicity, and introduce the following normalized endomorphisms
\bea
\tilde h_\epsilon={1\over {\rm sup}_X{\rm Tr}\,h_\epsilon}h_\epsilon.
\eea
Then $\|\tilde h_\epsilon\|_{L^\infty}\leq 1$. The key inequality that they satisfy is the following
\bea
\|\na\tilde h_\epsilon^\sigma\|_{L^2}^2\leq C,
\qquad \sigma\in [0,1],
\eea
where $C$ is independent of both $\epsilon$ and $\sigma$. Applying Rellich's Lemma and by a delicate diagonalization argument, one can then show that, again after replacing the $\epsilon$'s by a subsequence if necessary, there exists endomorphisms ${\mathcal H}_\sigma\in W^{2,1}(X,End\, E)$ for any $\sigma\in {\bf Q}\cap [0,1]$ such that
\bea
&&
\tilde h_\epsilon^\sigma\to {\mathcal H}_\sigma\ \ {\rm a.e.\ and\ in}\ L^2; \nonumber\\
&&
\na \tilde h_\epsilon^\sigma\to \na{\mathcal H}_\sigma\ \ {\rm weakly\ in}\ L^2,
\eea
as $\epsilon\to 0$. From this, one obtains the endomorphism $\pi\in W^{2,1}(X,End\,E)$ by 
\bea
\pi=I-{\rm lim}_{\sigma\to 0}{\mathcal H}_\sigma
\eea
which is not identically $0$, and satisfies the crucial properties
\bea
\label{UYpi}
&&\pi=\pi^*=\pi^2\nonumber\\
&&(1-\pi)\bar\pl\pi=0 \ \ a.e.
\eea

\medskip
$\bullet$ Remarkably, these properties of $\pi$ suffice to produce the desired sheaf. This is a consequence of a theorem on separate meromorphicity which is proved by Uhlenbeck and Yau in their paper \cite{UhYa}. For our discussion, it is convenient to adopt the formulation of this theorem by Shiffman \cite{Sh} and Popovici \cite{Po}:

\medskip
Let $E\to X$ be a holomorphic vector bundle over a compact K\"ahler manifold $X$, and $H$ a Hermitian metric on $E$. If $\pi\in W^{2,1}(X,End\,E)$ satisfies the properties listed in (\ref{UYpi}), then there exists a coherent sheaf ${\mathcal F}\subset {\mathcal O}(E)$ and an analytic subvariety $Z\subset X$of codimension $\geq 2$ such that

\bigskip

(a) $\pi_{\vert_{X\setminus Z}}\in C^\infty(X\setminus Z,End \,E)$;

(b) $\pi=\pi^*=\pi^2$ and
$(1-\pi)\bar\pl\pi=0\ \ {\rm on}\ X\setminus Z$;

(c) ${\mathcal F}_{\vert_{X\setminus Z}}=\pi_{\vert_{X\setminus Z}}(E)\to X\setminus Z$ is a holomorphic subbundle of $E_{\vert_{X\setminus Z}}$.

\bigskip
$\bullet$ Applying this separate meromorphicity theorem to the projection $\pi$ constructed in the previous steps from a sequence of endomorphisms $\tilde h_\epsilon$
gives a sheaf ${\mathcal F}$ that can be readily verified to be destabilizing. This completes the proof.

\bigskip
The Donaldson-Uhlenbeck-Yau theorem has now been extended to allow singularities, notably by Bando and Siu \cite{BS}, and very recently by Paun et al \cite{CGN}. This last work relied in an essential way on the mean-value inequality in (\ref{eqn:2.60}). Even more recently, the Bogomolov-Gieseker inequality has been extended by Guenancia and Paun \cite{GuPa} to reflexive ${\bf Q}$-sheaves in K\"ahler 3-folds with log-terminal singularities,
confirming a conjecture of Campana, H\"oring, and Peternell. However, to the best of the authors' knowledge, the Uhlenbeck-Yau method itself has not yet been applied to other equations, despite the fact that many continue to arise whose solvability has been conjectured to be equivalent to an algebraic stability condition. See notably the Hull-Strominger system and the dHYM equation discussed below. One can hope that this situation will change, perhaps with equations from symplectic geometry such as the Hitchin gradient flow equation which requires a short-time regularization reminiscent of the one used by Uhlenbeck and Yau, as discussed in \cite{FP}.

\bigskip

\subsection{\bf COMPLETE CALABI-YAU METRICS}
\hfil\break

Very early on after his solution of the Calabi conjecture, at the 1978 International Congress of Mathematicians in Helsinki, Yau \cite{Y80} had already called attention to the problem of finding complete Calabi-Yau metrics. This problem turns out to be quite hard, but there has been some breakthroughs recently, which suggests that some accelerating progress may be around the corner, with deep relations with other areas of the theory of partial differential equations:

\bigskip
{\bf (a) The Tian-Yau metric}

The first major result was due to Yau himself around 1990, in joint work with Gang Tian:

\bigskip
{\bf\small Theorem of Tian-Yau} \cite{TY89, TY90}: Let $X$ be a compact smooth Fano manifold. If $D$ is a smooth irreducible anti-canonical divisor on $X$, then the complement $X\setminus D$ admits a complete Calabi-Yau metric.

\bigskip
The Tian-Yau method is to find a candidate for the asymptotics of the desired Calabi-Yau metric near the divisor. Once such a candidate has been identified, they developed PDE methods which can then correct the asymptotic candidate to the desired Calabi-Yau metric. The Tian-Yau methods have since been refined by Hein \cite{H} and other authors, and give us some confidence that the most difficult step in the search of complete Calabi-Yau metrics is the first step of finding a suitable candidate for the asymptotics.

\bigskip
{\bf (b) The general problem of a divisor with normal crossings}

However, ideally, one would like to allow the divisor $D$ to have normal crossings. There was practically no progress on this problem until around 2020, where Tristan Collins and Yang Li were able to establish the existence of a complete Calabi-Yau metric under the following conditions:

\bigskip
{\bf \small Theorem of Collins-Li} \cite{CL1}: Let $X$ be a smooth Fano manifold of dimension $n\geq 3$. Assume that its anti-canonical divisor is of the form $(d_1+d_2)L$, where $L$ is a positive line bundle, and $d_1$ and $d_2$ are two positive integers. Let $D_1$ and $D_2$ be two transversally intersecting smooth divisors in the linear systems associated to $d_1L$ and $d_2L$ respectively. Then $X\setminus(D_1\cup D_2)$ admits a complete Calabi-Yau metric.

\medskip
An important innovation in \cite{CL1} is to draw on the Calabi ansatz and on torus invariant dimensional reductions of Calabi-Yau metrics to produce an ODE for a candidate asymptotics near the divisors. Not only did this turn out to lead indeed to the required asymptotics, but it also suggested that, for more component divisors and by an inductive procedure on the number of components, the generalization of this ODE will be some {\it real} Monge-Amp\`ere equation.

\bigskip
{\bf (b) Relation with free-boundary problems and optimal transport}

\medskip
Motivated by these heuristics, the following boundary value problem was proposed by Collins, Tong, and Yau \cite{CoToYa}:

\bigskip
Let $P\subset{\bf R}^n$ be an open convex set containing the origin and $k$ be a positive integer $\leq n$. Find a convex function $v$ satisfying
\bea
\label{MAreal}
{\rm det}\,D^2v&=&v^{-(n+2)}(-v^*)^{-k}\ \ {\rm on}\ P\nonumber\\
v^*&=&0 \ \ {\rm on}\ \pl P.
\eea
where $v^*$ is the composition of the Legendre transform of $v$ with the gradient map $\na v$, 
\bea
v^*(y)=\<y,\na v(y)\>-v(y).
\eea

\bigskip
Similar equations with the more familiar Dirichlet conditions have been studied in \cite{ToYa}. A novelty in the above problem is rather the boundary condition on $v^*$, which is equivalent to, as explained in \cite{CoToYa}, a free-boundary problem on the Legendre transform $u(x)={\rm sup}_P(\<x,y\>-v(y))$,
\bea
{\rm det}\,D^2u&=&(u^*)^{n+2}{\rm max}(-u,0)^k\ \ {\rm on}\ \ {\bf R}^n
\nonumber\\
\na u({\bf R}^n)&=& \bar P.
\eea
We have then:

\bigskip
{\bf \small Theorem of Collins-Tong-Yau} \cite{CoToYa} : Let $P\subset{\bf R}^n$ be an open convex set containing the origin, and $1\leq k\leq n$ be an integer. Then the above problem has a unique solution $v\in C^\infty(P)\cap C^{1,\alpha}(\bar P)$ for some $\alpha>0$.

\bigskip

Note that the solution cannot be in $C^2$ up to the boundary, as the equation is singular there.

\bigskip
Remarkably, these developments have inspired Collins and Tong \cite{CoTo} to a new method to study the regularity of optimal transport maps, using a monotonicity formula, which allowed them on one hand, to extend a regularity theorem of Savin-Yu \cite{SaYu} from two to all dimensions, and on the other hand, to extend a theorem of Chen-Liu-Wang \cite{CLW} from $C^{1,1}$ convex domains to $C^{1,\beta}$ domains. These improved versions are now sharp.

\bigskip
As of the writing of the present survey paper, these are very recent developments still in flux. Nevertheless, it is clear that a problem formulated some 50 years ago by Yau has brought new relations between two seemingly distant subjects, namely complex geometry and optimal transport, enriching immensely both subjects.

\bigskip

\subsection{\bf HULL-STROMINGER EQUATIONS}
\hfill\break

The Ricci-flat equation and the Hermite-Einstein equation are K\"ahler versions of the equations for gravitation and gauge theories, which are field theories describing individual fundamental forces of nature. But it has been a discovery in 1985 with very wide ramifications for both mathematics and physics, by Candelas, Horowitz, Strominger, and Witten \cite{CHSW}, that the Ricci-flat condition arises also from unified string theories, this time from a different type of requirement, namely supersymmetry. For our purposes, it suffices to know that supersymmetry is a symmetry generated by a spinor field. Now 
unified string theories are supersymmetric theories of extended objects which take place in a 10-dimensional Lorentz space-time. To make contact with the 4-dimensional Lorentz space-time $M^{1,3}$ of our common day experience, the 10-dimensional space-time is taken to be of the form $M^{1,3}\times X$, where $X$ is a very small $6$-dimensional 
manifold, a process known as compactification.
For phenomenological reasons, one requires that supersymmetry be unbroken in this compactification. This implies that the spinor field generating supersymmetry is covariantly constant with respect to a Levi-Civita connection suitably modified by torsion terms encoded in a $3$-form known as flux.
As discovered in \cite{CHSW}, one way of insuring this is to set the flux to $0$ and require that the internal space $X$ be a Calabi-Yau $3$-fold.

\bigskip
Very shortly after the appearance of \cite{CHSW}, a generalization to non-K\"ahler manifolds was proposed independently by Hull \cite{Hu} and Strominger \cite{St}. Let $X$ be a compact complex $3$-fold, equipped with a nowhere vanishing holomorphic $(3,0)$-form. Let $E\to X$ be a holomorphic vector bundle. Then the Hull-Strominger system is the following system of equations for a Hermitian metric $\o$ on $X$ and a Hermitian metric $h$ on $E$,
\bea
\label{HuSt}
&&
\o^2\wedge F_h=0\nonumber\\
&&i\ddb\, \o={\alpha'\over 4}{\rm Tr}(Rm\wedge Rm -F_h\wedge F_h)\\
&& d(\|\Omega\|\o^2)=0,\nonumber
\eea
where $Rm$ and $F_h$ denotes the curvature forms of the Chern unitary connection of $\o$ and $h$ respectively, and $\alpha'$ is a constant\footnote{Other unitary connections for $\o$ have also been considered in the literature.}. Note that the third equation is a conformal version of the balanced condition of Michelsohn \cite{Mi}. Clearly the Hull-Strominger system requires the consistency conditions in Bott-Chern cohomology,
\bea
&& c_1(E)=0\in H^{1,1}_{BC}(X,{\bf R})\nonumber\\
&& c_2(E)=c_2(X)\in H^{2,2}_{BC}(X,{\bf R}).\nonumber
\eea
Calabi-Yau manifolds are a special solution of this system, corresponding to taking $E=T^{1,0}(X)$, $\o=h$, and $\o$ K\"ahler. Indeed, the first equation is satisfied if $\o=h$ and $\o$ is Ricci-flat. The Ricci-flatness of $\o$ implies that $\|\Omega\|$ is constant, which implies the third equation when combined with the K\"ahler property of $\o$. Finally, the second equation is also trivially satisfied, as the right hand side vanishes for $\o=h$ and the left hand side vanishes for K\"ahler $\o$. 
Thus the Hull-Strominger system can be viewed as a generalization of Calabi-Yau manifolds which allows for metrics with non-zero torsion.

\bigskip
The Hull-Strominger system appeared very complicated and did not seem to attract much attention at the start. However, the situation changed drastically after the discovery by J.X. Fu and Yau around 2006 \cite{FY1, FY2} of a non-K\"ahler solution by PDE methods. This solution is built on an adaptation by Goldstein and Prokushkin of a classical construction of Calabi and Gray, and is a torus fibration over a $K3$ surface (with the holomorphic bundle $E$ taken to be flat). In this case, the quadratic terms in the curvature simplified considerably, and the equation reduces to a Monge-Amp\`ere equation, which they managed to solve. Since that time, many more solutions have been found, see e.g. \cite{Fe, FeYa, AGF, FGV} and references therein.
This revealed that the Hull-Strominger system has a very rich mathematical structure, and its interest lies beyond the original goal of enlarging the space of possible supersymmetric compactifications for the heterotic string. Too many research avenues originating from the Hull-Strominger systems have been explored since, and we can only mention a few.

\bigskip
A first avenue has been suggested by Yau, regarding the well-known ``Reid's fantasy". This is a proposal by Reid \cite{Re} that moduli spaces of Calabi-Yau 3-folds with different Hodge numbers could be connected by conifold transitions, which are topology changing processes discovered by Clemens \cite{Cl} and Friedman \cite{Fr}, and that the resulting web of Calabi-Yau spaces is connected. The question raised by Yau is whether this proposal can be implemented at the level of canonical metrics. Since conifold transitions may not preserve the K\"ahler property, we would need a notion of canonical metric defined by an equation admitting non-K\"ahler solutions. A candidate for such a notion could be the solution of Hull-Strominger systems. This avenue is clearly very challenging, but it is being vigorously pursued \cite{CGPY}. See in particular the lectures of Collins \cite{Co}.

\bigskip
Another very intriguing question is to determine necessary and sufficient conditions for the solvability of Hull-Strominger systems, which are still not known at the moment. Some suitable stability conditions may be the right ones, but this is still an open question. This is being explored in particular by M. Garcia-Fernandez and his co-workers \cite{GFM,GFM1, GRST}, drawing inspiration from both the theory of vertex algebras and generalized geometry.

\bigskip
From the more analytic viewpoint, the Fu-Yau solution of the Hull-Strominger system suggests some particular Hessian equations, now known as Fu-Yau equations, which may be interesting in their own right \cite{PPZ3}. Solutions for these equations have been obtained in certain regimes, but the general solution is again very far from being understood.

\bigskip
Finally we mention a broad issue with analysis on non-K\"ahler manifolds and which will occur repeatedly below. It has to do with the fact that, at least in equations arising from unified string theories, the K\"ahler condition does not disappear altogether, but is rather replaced by a weaker cohomological condition. The conformally balanced condition in Hull-Strominger systems is a prime example. The issue is how  to implement such conditions in the absence of a $\pl\bar\pl$-Lemma. It has been advocated in \cite{PPZ1,PPZ2} that the most effective way may be by a geometric flow which preserves the desired condition, and start from an initial data satisfying it. See also \cite{Bryant, BeVe} and \cite{GSM}. These flows turn out to be natural generalizations of the K\"ahler-Ricci flow to the non-K\"ahler setting. Thus Hull-Strominger systems also provide a useful laboratory for developing PDE methods for geometric flows.

\bigskip
\subsection{MIRROR SYMMETRY AND THE STROMINGER-YAU-ZASLOW CONJECTURE}

\hfill\break

\smallskip
Mirror symmetry burst upon the scene in 1989, and has been a source of challenges and inspiration for practically all branches of mathematics ever since. Very roughly speaking, a Calabi-Yau manifold and its mirror are manifolds defining the same quantum field theory. A precise mathematical version of this characterization is obviously hard to formulate. A celebrated version was proposed in 1996 by Strominger-Yau-Zaslow \cite{SYZ} and it has been a powerful engine for many developments in both complex geometry and symplectic geometry. We reproduce here the version given in \cite{CoLin}:

\bigskip
{\bf \small The SYZ Conjecture}: Let $(X^\vee,\o^\vee)$ be a Calabi-Yau manifold, and let ${\mathcal M}_{cplx}^\vee$ be the moduli space of complex structures on $X^\vee$. If $J^\vee$ is a complex structure on $X^\vee$ sufficiently close to a large complex structure limit, then

\medskip
(a) $(X^\vee, J^\vee,\o^\vee)$ admits a special Lagrangian torus fibration $\pi^\vee:X^\vee\to B^\vee$ onto a base $B^\vee$ equipped with an integral affine structure;

\medskip
(b) There is another Calabi-Yau manifold $(X,J,\o)$ with a special Lagrangian fibration $\pi:X\to B$ over a base $B$ equipped with an integral affine structure;

\medskip
(c) Let ${\mathcal M}_{Kah}$ be the complexified K\"ahler moduli space of $X$. Then there is a mirror map $q: {\mathcal M}_{cplx}^\vee\to {\mathcal M}_{Kah}$ which is a local diffeomorphism and which satisfies $q(J^\vee)=\o$;

\medskip
(d) There is an isomorphism $\varphi: B^\vee\to B$ exchanging the complex and affine symplectic structures, and such that the Riemannian volumes of the Lagrangian torus fibers over $b^\vee\in B^\vee, \varphi(b^\vee)\in B$ are inverses of each other.

\bigskip
We note the correspondence between complex structures and symplectic structures, which is one of the characterizing features of mirror symmetry. The correspondence between the torus fibers is usually known as $T$-duality.
An idea of how difficult the SYZ conjecture is can be gathered from the fact that, in general, it is already quite difficult to find a single special Lagrangian torus inside a Calabi-Yau 3-fold, let alone a whole fibration. Nevertheless, there has been in the last few years two remarkable advances: on one hand, the SYZ conjecture, as formulated above, was proved by T. Collins, A. Jacob  and J. Lin \cite{JaCoLi1, JaCoLi2} for certain complete, non-compact, 2-folds, of the form $X=Y\setminus D$, where $Y$ is a Del Pezzo surface, and $D\in |-K_Y|$ is a smooth elliptic curve. It may be noteworthy that the Tian-Yau metric described earlier played a major role here; on the other hand, a weaker statement than the full SYZ conjecture was established by Yang Li \cite{Li} for the Fermat family
\bea
X_s=
\{Z_0Z_1\cdots Z_n+e^{-s}\sum_{j=0}^{n+1}Z_j^{n+2}=0\},
\qquad s>>1.
\eea
He proves that there exists a subsequence $X_s$ with $s\to\infty$ which admits a Lagrangian torus fibration over a generic region $U_s\subset X_s$, with the ratio of volume of $U_s$ to volume of $X_s$ tending to $1$ as $s\to\infty$. 

We refer to \cite{JaCoLi1,JaCoLi2,Li, CoLin} for fuller statements. It is tempting to believe that, even with the challenging nature of the SYZ conjecture, progress may soon be accelerating.

\bigskip

\subsection{SYMPLECTIC COHOMOLOGY AND STRING EQUATIONS MODELS}

\hfill\break

A notion of symplectic Hodge theory had been introduced a long time ago by Ehreshmann and Liebermann \cite{EhLi} and Brylinski \cite{Br}, but it suffers from several drawbacks, notably the lack of existence and uniqueness within a cohomology class of the corresponding analogue of harmonic forms. This problem was addressed by Li-Sheng Tseng
and Yau \cite{LiYa1, LiYa2}, who introduced new symplectic cohomologies  and also showed how their theories can be applied to string theory \cite{LiYa3}. For the more analytic issues discussed in the present paper, a particularly valuable contribution of their work \cite{LiYa3} is the fact that they wrote down explicitly the equations for the Type IIA 
and the Type IIB string theories with O6/D6 and O5/D5 brane sources respectively. These equations suggested some new flows, the study of which revealed some important structures, in particular for symplectic geometry.

\bigskip
Here we would like to mention the following so-called Type IIA flow, which is a flow of closed, primitive, $3$-forms $\varphi$ on a compact $6$-dimensional symplectic manifold $(X,\o)$ introduced by Fei, Phong, Picard, and Zhang in \cite{FPPZ},
\bea
\dot\varphi&=&d\Lambda d(|\varphi|^2\star\varphi)\nonumber\\
\varphi(0)&=&\varphi_0,
\eea
the stationary points of which satisfy the equation without source written by Tseng and Yau \cite{LiYa3}.
We recall that, on a general $6$-dimensional manifold, Hitchin \cite{Hi} had shown how, by a pure pointwise and algebraic construction, a $3$-form $\varphi$ would give rise to an almost-complex structure $J_\varphi$. On a symplectic manifold $(X,\o)$, if we require $\varphi$ to be primitive, that is, if $\Lambda\varphi=0$ where $\Lambda:\wedge^n(X)\to\wedge^{n-2}(X)$ is the Hodge operator of contracting with $\o$, we can then obtain a Hermitian form
\bea
g_\varphi(U,V)=\o(U,J_\varphi V)
\eea
which is a Hermitian metric under the open condition that it is positive. Thus each generic $\varphi$ gives rise to an almost-K\"ahler manifold $(X,\o,J_\varphi)$. When $\varphi$ is closed, we refer to $(X,\o,\varphi,J_\varphi)$ as a Type IIA structure. It is one of the important results of \cite{FPPZ} that a Type IIA structure is an almost-Hermitian manifold with SU(3) holonomy, with the key distinction that it is with respect to the projected Levi-Civita connection, and not the Levi-Civita connection. This structure is crucial for the existence and Shi-type estimates for the above Type IIA flow. In particular, it is needed for the square of the norm of the Nijenhuis tensor of $J_\varphi$ to obey a parabolic flow, which is an important property not shared by many flows of almost-complex structures.

\bigskip
Our interest in discussing the Type IIA flow, besides its arising from another recent direction opened up by Yau, lies in its likely relations with some issues that we had discussed earlier, namely both the Uhlenbeck-Yau method and free boundary problems. It is indeed natural to ask whether the short-time existence of certain flows in symplectic geometry \cite{FP}, such as the Hitchin gradient flow which is still an open problem, can be obtained through a regularization scheme similar to the one used by Uhlenbeck-Yau in their solution of the Hermitian-Einstein equation. Furthermore, in the above discussion of the Type IIA flow, we had not incorporated as yet any source. We expect that sources will arise from a free boundary problem as discussed in the problem of Tian-Yau metrics. Clearly there is a lot that needs to be investigated.

\bigskip

\subsection{\bf THE dHYM EQUATION}

\hfill\break

\smallskip
The deformed Hermitian-Yang-Mills (or dHYM) equation is the following equation. Let $(X,\o)$ be an $n$-dimensional compact, connected, K\"ahler manifold of dimension $n$, and let $\chi$ be a closed real $(1,1)$-form on $X$. The question is to determine whether there exists a smooth $(1,1)$-form $\chi_u=\chi+i\ddb u\in [\chi]$
satisfying the equation
\bea
\label{dHYM}
\sum_{j=1}^n{\rm arctan}\,\lambda_j=\hat\theta .
\eea
Here $\{\lambda_j\}_{j=1}^n=\{\lambda([\chi_u])\}_{j=1}^n$ are the eigenvalues of the endomorphism $h^j{}_k
=\o^{j\bar m}(\chi_u)_{\bar mk}$, and $\hat\theta$ is a topological constant. This equation stands at the crossroads of the most active areas in geometry and physics: it was motivated by mirror symmetry, and proposed independently by Marino, Minasian, Moore, and Strominger \cite{MMMS} and by Leung, Yau, and Zaslow \cite{LYZ}. It appeared first in the mathematical literature in the 2017 paper by A. Jacob and Yau \cite{JY} who solved it in dimension $n=2$, where it can be reduced to a Monge-Amp\`ere equation. It is related by mirror symmetry to the existence of special Lagrangian manifolds in the mirror. It is also the natural K\"ahler analogue of the special Lagrangian equation introduced by Harvey and Lawson \cite{HL2} in the symplectic setting. The equation is said to be supercritical if $\hat\theta\in (0,\pi)$. The well-known $J$-equation introduced earlier by X.X. Chen \cite{ChX} and Donaldson \cite{Do2} can be viewed as the small radius limit of the dHYM equation.

\bigskip

The following theorem gives criteria for the existence of solutions in dimensions $n\geq 3$ in terms of subsolutions:
\bigskip

{\bf \small Theorem of Collins-Jacob-Yau}
\cite{CoJaYa}: Let $(X,\o,\chi)$ be as above, and $\theta_0\in (0,\pi)$. Assume that there exists a function $\underline u$ satisfying the following two conditions:

\medskip

(a) $\underline u$ is a subsolution of the dHYM equation in the sense of Guan \cite{Gu} and Szekelyhidi \cite{Sze}, that is
\bea
{\rm sup}_X{\rm sup}_{1\leq j\leq n}\sum_{k\not=j}{\rm arctan}\,\lambda_k(\chi_{\underline u})<\theta_0;
\eea

(b) $\chi({\underline u})$ also satisfies the inequality
\bea
{\rm sup}_X\,\sum_{j=1}^n {\rm arctan}\,\lambda_j(\chi_{\underline u})<\theta_0.
\eea
Then the dHYM equation admits a unique smooth solution.

\bigskip
Collins-Jacob-Yau also conjectured that the second condition (b) can be removed. This was achieved by Pingali \cite{Pi} for $n=3$ and by Lin \cite{Li} for $n=4$. A flow based proof of the above theorem was subsequently given by Fu, Yau, and Zhang \cite{FYZ}.
There the notion of subsolution is the parabolic version given by Phong-To \cite{PhTo}.

\bigskip
A major question remains, which is whether the subsolution condition can be replaced by a suitable algebraic stability condition in the sense of GIT. For the $J$ equation, which is the small radius limit of the dHYM equation, the solvability in terms of a subsolution was a result of Song and Weinkove \cite{SoWe} from 2008, an algebraic stability characterization was conjectured by Lejmi and Szekelyhidi \cite{LS} in 2015,
and proved only rather recently by Gao Chen \cite{ChGa} and Song \cite{S20} around 2020. The algebraic stability condition for the dHYM appears to be a formidable problem, which is being explored by many authors. Surveys can be found in \cite{CoSh,CLSY}. Here we mention briefly only a few results.

\bigskip
In their original paper \cite{CoJaYa}, Collins, Jacob, and Yau had investigated in detail algebraic conditions for the solvability of the dHYM equation when $X$ is a surface. In particular, they showed that the solvability condition is equivalent to the following ``twisted ampleness" condition: let $X$ be a K\"ahler surface, and $L\to X$ a holomorphic line bundle. For any curve $C\subset X$, define the following complex numbers,
\bea
Z_X(L)=-\int_X e^{-\sqrt{-1}\o}ch(L),
\qquad
Z_C(L)=-\int_C e^{-\sqrt{-1}\o}ch(L).
\eea
Assume that $\Im \,Z_X(L)>0$. Then the dHYM equation admits a solution if and only if for any curve $C\subset X$, we have
\bea
\label{twistedampleness}
\Im({Z_C(L)\over Z_X(L)})>0.
\eea
In $n$ dimensions, a first result was obtained by Jacob and Sheu \cite{JaSh} in the case of the blow-up ${\bf Bl}_p{\bf P}^n$ at one point.
We should also mention that, 
motivated by the results of Chen and Song for the $J$-flow, a Nakai-Moishezon criterion for the dHYM equation was also introduced by Chu, Lee, and Takahashi \cite{ChLeTa}.

\bigskip

However, a very important and intriguing question is whether one can characterize the solvability of the dHYM equation by a Bridgeland stability condition. These stability conditions had been introduced by Bridgeland \cite{Br1} based on a proposal in physics of M. Douglas \cite{Dou}.
Since the dHYM equation is heuristically dual to the problem of existence of special Lagrangians, one may look for guidance to the special Lagrangian problem. Already in the early 2000's, it had been proposed by Thomas and Yau \cite{ThYa} that the existence of special Lagrangians and the convergence of the Lagrangian mean curvature flow should be equivalent to an algebraic stability condition. The precise notion of stability needed was left open. More recently, it had been conjectured by Joyce \cite{Jo} that a Bridgeland stability condition would be the right one for the special Lagrangian problem. Thus it is natural to expect Bridgeland stability conditions to be at least related in some way to the existence of solutions of the dHYM equation. Indeed, as already noted by \cite{CoJaYa}, the number $Z_X(L)$ is reminiscent of the notion of central charge in Bridgeland stability conditions.
The problem raised by Collins and Yau of determining precisely the relations between the solvability of dHYM, Bridgeland stability, and other stability conditions, appears rather subtle.
For example, using the work of Arcara and Miles \cite{ArMi}, it was shown by Collins and Shi \cite{CoSh} that a line bundle $L\to {\bf Bl}_p{\bf P}^2$ admitting a metric solving the dHYM equation is Bridgeland stable, but not conversely. Alternatively, for the last few years, Collins and Yau \cite{CoYa} have embarked on an ambitious program to address the dHYM equation through an infinite-dimensional version of GIT theory. The ultimate  stability condition should arise from this program.

\setcounter{equation}{0}


\begin{thebibliography}{99}

{\footnotesize



\bibitem{AGF} B. Andreas and M. Garcia-Fernandez,
{\em Solutions to the Strominger system via stable bundles over Calabi-Yau threefolds},
Commun. Math. Phys. 315 (2012) 153-168.

\bibitem{ArMi} D. Arcara  and  E. Miles, E. {\em Bridgeland stability of line bundles on surfaces}, J. Pure App. Alg., 220, (2016), no. 4, 1655-1677.

\bibitem{BS}  S. Bando and Y.T. Siu,   {\em Stable sheaves and Einstein-Hermitian metrics}, in Geometry and Analysis on Complex Manifolds, World Sci. Publ., River Edge, NJ, 1994, 39-50.

\bibitem{BeVe} L. Bedulli and L. Vezzoni,
{\em On the stability of the anomaly flow}, Math. Res. Lett. 29 (2022) 323-338.

\bibitem{Be} F. Bei, {\em On the Laplace-Beltrami operator on compact complex spaces},  Trans. Amer. Math. Soc. 372 (2019), no. 12, 8477--8505.

\bibitem{Bl} Z. Blocki, {\em On the uniform estimate in the Calabi-Yau theorem II}, Science China Math. 54 (2011) 1375-1377.

\bibitem{Br1} T. Bridgeland, {\em Stability conditions on triangulated categories}, Ann. of Math. (2) 166 (2007), no. 2, 317-345.

\bibitem{Bryant} R. Bryant and F. Xu,
{\em Laplacian flow for closed $G2$ structures: short time behavior},
arXiv: 1101.2004.

\bibitem{Br} J.L. Brylinski, J.-L. {\em A differential complex for Poisson manifolds,}
J. Differential Geom. 28(1): 93-114 (1988).


\bibitem{CHSW} P. Candelas, G. Horowitz, A. Strominger, and E. Witten,
{\rm Vacuum configurations for superstrings}, Nucl. Phys. B 258 (1985) 46-74.

\bibitem{CGN} J. Cao, P. Graf, P.   Naumann, M. Paun, T. Peternell, and  X. Wu,  {\em Hermitian--Einstein metrics in singular settings}, arXiv:2303.08773.

\bibitem{ChGa} G. Chen, {\em The $J$-equation and the supercritical deformed Hermitian-Yang-Mills equation}, Invent. math. 225, 529-602 (2021).

\bibitem{ChX} X.X. Chen,  {\em A new parabolic flow in K\"ahler manifolds}, Comm. Anal. Geom. 12 (2004), no. 4, 837 - 852.

\bibitem{CC}X.X. Chen and J.R. Cheng, {\em On the constant scalar curvature K\"ahler metrics I - a priori estimates}, J. Amer. Math. Soc. (2021) DOI: https://doi.org/10.1090/jams/967.

\bibitem{CLW} S. Chen, J. Lu, and X.J. Wang,
{\em Global regularity for the Monge-Amp\`ere equation with natural boundary condition},
Ann. of Math. 194 (2021) no. 3, 745-793.

\bibitem{CL} S.Y. Cheng and P. Li,   {\em Heat kernel estimates and lower bound of eigenvalues}, Comment. Math. Helvetici 56 (1981) 327-338.


\bibitem{ChLeTa} J. Chu, M. Lee, and R. Takahashi, {\em A Nakai-Moishezon type criterion for supercritical deformed Hermitian-Yang-Mills equation,} J. Differ. Geom. 126.2 (2024): 583-632.

\bibitem{Cl} H. Clemens,
{\em Double solids}, Adv. Math. 47 (1983) no. 2, 107-230.

\bibitem{Co} T. Collins, {\em Introduction to conifold transitions}, arXiv:2509.01002.

\bibitem{CGPY} T. Collins, S. Gukov, S. Picard, and  S.T. Yau, {\em Special Lagrangian cycles and Calabi-Yau transitions}, Comm. Math. Phys., 401 (2023), no. 1, 769 - 802.

\bibitem{JaCoLi1} T. Collins, A. Jacob, and Y. Lin, {\em Special Lagrangian submanifolds of log Calabi-Yau manifolds}, Duke Math. J., 170 (7), 1291-1375.

\bibitem{JaCoLi2} T. Collins, A. Jacob, and Y. Lin, {\em The SYZ mirror symmetry conjecture for del Pezzo surfaces and rational elliptic surfaces}, arXiv:2012.05416.

\bibitem{CoJaYa} T. Collins, A. Jacob, and S.T. Yau,
{\em (1,1)-forms with specified Lagrangian phase: a priori estimates and algebraic obstructions}, Cambridge J. Mathematics Vol. 8 no. 2 (2020) 407-452.

\bibitem{CL1}  T. Collins and Y. Li, {\em Complete Calabi-Yau metrics in the complement of two divisors}, arXiv:2203.10656.

\bibitem{CoLin} T. Collins and Y.S. Lin,
{\em Recent progress on SYZ mirror symmetry for some non-compact Calabi-Yau surfaces}, arXiv: 2208.14485.



\bibitem{CLSY} T. Collins, J. Lo, Y. Shi, and S.T. Yau,
{\em Stability for line bundles and deformed Hermitian Yang-Mills equation on some elliptic surfaces}, arXiv: 2306.05620.


\bibitem{CoSh} T. Collins and Y. Shi, {\em Stability and the deformed Hermitian-Yang-Mills equation}, Surveys in Differential Geometry, 24(1) (2021), 1-38.  

\bibitem{CoTo} T. Collins and  F. Tong, {\em Boundary regularity of optimal transport maps on convex domains}, arXiv:2507.05395.

\bibitem{CoYa} T. Collins and S.T. Yau, {\em Moment maps, nonlinear PDE, and stability in mirror symmetry I: Geodesics},  Annals of PDE (1) 11 (2021) 68 pp.

\bibitem{CoToYa} T. Collins, F. Tong, and S.T. Yau, {\em A free boundary Monge-Amp\`ere equation and applications to complete Calabi-Yau metrics}, arXiv:2402.10111.

\bibitem{CMM} D. Coman, X. Ma, and G. Marinescu, {\em Equidistribution for sequences of line bundles on normal K\"ahler spaces}, Geometry \& Topology 21 (2017) 923 - 962.

\bibitem{DP} J.P. Demailly and N.Pali, {\em Degenerate complex Monge-Amp\`ere equations over compact K\"ahler manifolds}, Internat. J. Math. 21 (2010) no. 3, 357-405.

\bibitem{DK} S. Dinew and S. Kolodziej, {\em A priori estimates for complex Hessian equations}, Anal. PDE 7 no 1 (2013) 227-244.

\bibitem{Do} S. Donaldson,  {\em Anti-self-dual Yang-Mills connections over complex algebraic surfaces and stable vector bundles}, Proc. LMS 50 (1985) 1-26

\bibitem{Do2} S. Donaldson,  {\em Moment maps and diffeomorphisms}, Asian J. Math., 3 (1999), 1-16.

\bibitem{DS} S. Donaldson and S. Sun, {\em Gromov-Hausdorff limits of K\"ahler manifolds and algebraic geometry}, Acta Math. 213 (2014), no. 1, 63-106.

\bibitem{Dou} M. Douglas, {\em Dirichlet branes, homological mirror symmetry, and stability}, Proceedings of the International Congress of Mathematicians, Vol. III (Beijing, 2002), 395-408, Higher Ed. Press, Beijing, 2002.

\bibitem{EhLi} A. Ehresmann, and D. Libermann, {\em Sur les structures presque hermitiennes isotropes},
Colloque de G\'eom\'etrie Diff\'erentielle Globale (Bruxelles, 1958),
Centre Belge de Recherches Math\'ematiques,
Gauthier-Villars, Paris, 1959, pp. 59 - 77.

\bibitem{EGZ} P. Eyssidieux, V. Guedj, and A. Zeriahi, {\em Singular K\"ahler-Einstein metrics}, J. Amer. Math. Soc. 22 (2009), 607-639.

\bibitem{Fe} T. Fei, {\em A construction of non-K\"ahler Calabi-Yau manifolds and new solutions to the Strominger system}, Adv. Math. 302 (2016) 529-550.

\bibitem{FP} T. Fei and D.H. Phong, {\em Symplectic geometric flows}, Pure Appl. Math. Quart. 19 (2023) no. 4, 1853-1871.

\bibitem{FPPZ} T. Fei, D.H. Phong, S. Picard, and X. Zhang, X. {\em Geometric Flows for the Type IIA String}, Cambridge J. Math. Vol. 9 no. 3 (2021) 693-807,
arXiv:2011.03662.

\bibitem{FeYa} T. Fei and S.T. Yau,
{\em Invariant solutions to the Strominger system on complex Lie groups and their quotients}, Commun. Math. Phys. 338 (2015) no. 3, 1183-1195.

\bibitem{FGV} A. Fino, G. Grantcharov, and L. Vezzoni,
{\em Solutions to the Hull-Strominger system with torus symmetry},
Commun. Math. Phys. 388 (2021) 947-967.

\bibitem{Fr} R. Friedman, {\em Simultaneous resolution of threefold double points},
Math. Ann. 274 (1986) no. 4, 671-689.

\bibitem{FGS} X. Fu, B. Guo, and J. Song, {\em Geometric estimates for complex Monge-Amp\`ere equations},  J. Reine Angew. Math. 765 (2020), 69-99.

\bibitem{FYZ} J.X. Fu, S.T. Yau, and D. Zhang,
{\em Introduction to a deformed Hermitian Yang-Mills flow},
Surveys in Diff. Geom. XXVI (2024) 157-168.

\bibitem{FY1} J.X. Fu and S.T. Yau, {\em The theory of superstring with flux on non-K\"ahler manifolds and the complex Monge-Amp\`ere equation} J. Differential Geom. 78 (2008), no. 3, 369 - 428.

\bibitem{FY2} J.X. Fu and S.T. Yau, S.-T. {\em A Monge-Amp\`ere type equation motivated by string theory}, Comm. Anal. Geom. 15 (2007) 29 - 76.

\bibitem{GFM1} M. Garcia-Fernandez and G. Gonzalez Molina,
{\em Harmonic metrics for the Hull-Strominger system and stability},
Inter. J. of Math. Vol 35 no 9 (2024) 2441008, 30 pp.

\bibitem{GFM} M. Garcia-Fernandez  and  G. Gonzalez Molina, {\em Futaki invariants and Yau’s conjecture on the Hull-Strominger system}, arXiv:2303.05274.

\bibitem{GRST} M. Garcia-Fernandez, R. Rubio, C. Shahbazi, and C. Tipler,
{\em Canonical metrics on holomorphic Courant algebroids}, Proc. London Math. Soc.
125 no. 3 (2022) 329-367.



\bibitem{GSM} M. Garcia-Fernandez, J. Jordan, and J. Streets, 
{\em Non-K\"ahler Calabi-Yau geometry }, J. Math. Pure Appl. 177 (2023) 329-367.

\bibitem{Gu} B. Guan, {\em Second order estimates and regularity for fully nonlinear elliptic equations on Riemannian manifolds}, Duke Math. J. 163 (2014), 1491-1524.

\bibitem{GT} V. Guedj and T.D. To, {\em K\"ahler families of Green's functions}, arXiv:2405.17232.

\bibitem{GuPa} H. Guenancia and M. Paun,
{\em Bogomolov-Gieseker inequality for threefolds with log terminal singularities},
arXiv: 2405.10003.


\bibitem{GP} B. Guo and D.H. Phong, {\em Auxiliary Monge-Amp\`ere equations in geometric analysis}, Notices of the ICCM, Volume 11 (2023) Number 1, 98-135.

\bibitem{GPa} B. Guo and D.H. Phong, {\em On $L^\infty$ estimates for fully nonlinear partial differential equations}, Ann. of Math. 200 (2024), no. 1, 365-398.

\bibitem{GPb} B. Guo and D.H. Phong,
{\em Uniform entropy and energy bounds for fully nonlinear equations},
Comm. Anal. Geom. 32 no. 8 (2024) 2305-2325.


\bibitem{GPSS23}  B. Guo, D.H. Phong, J.  Song,  and J. Sturm,  {\em Sobolev inequalities on K\"ahler spaces}, preprint, (2023), arXiv:2311.00221.


\bibitem{GPSS} B. Guo, D.H.  Phong, J.  Song,  and J. Sturm, {\em Diameter estimates in K\"ahler geometry}, Comm. Pure Appl. Math.,  Volume 77, Issue 8 (2024), 3520-3556.

\bibitem{GPSS24} B. Guo, D.H.  Phong, J. Song, and J. Sturm, {\em Diameter estimates in K\"ahler geometry II: removing the small degeneracy assumption},  Math. Z. 308, 43 (2024) arXiv:2405.18280.

\bibitem{GPS}  B. Guo, D.H. Phong, and J. Sturm, {\em Green's functions and complex Monge-Amp\`ere equations}, J. Differential Geom. Vol. 127, No. 3 (2024) 1083-1119.


\bibitem{GPT} B. Guo, D.H.  Phong, and F. Tong, {\em On $L^\infty$ estimates for complex Monge-Amp\`ere equations},  Ann. of Math. (2) 198 (2023), no.1, 393-418. 

\bibitem{HL} F.R. Harvey and H.B. Lawson, {\em Determinantal majorization and the work of Guo-Phong-Tong and Abja-Olive},  Calc. Var.  Partial Diff. Equations 62, 153 (2023).

\bibitem{HL1} F.R. Harvey and H.B. Lawson,
{\em A definitive majorization result for nonlinear operators},
Duke Math. J. 174 (13) (2025) 2749-2763.

\bibitem{HL2} F.R. Harvey and H.B. Lawson,  {\em Calibrated geometries}, Acta Math., 148 (1982), 47-157.

\bibitem{H} H.J. Hein, {\em Gravitational instantons from rational elliptic surfaces},
J. Amer. Math. Soc. 25 (2012) 355-393.

\bibitem{Hi} N. Hitchin,   {\em The geometry of three-forms in six dimensions}, J. Differential Geom. 55 (2000), no. 3, 547-576.

\bibitem{Ho}  L. H\"ormander, {\em An introduction to complex analysis in several variables}. Van Nostrand, Princeton, NJ, 1973.

\bibitem{Hu} C. Hull, {\em Compactifications of the heterotic superstring}, Phys. Lett. B 1978 (1986), no. 4, 357 - 364.

\bibitem{JaSh} A. Jacob and N. Sheu, {\em The deformed Hermitian-Yang-Mills equation on the blowup of $\mathbb {P}^n$}, arXiv:2009.00651.

\bibitem{JY} A. Jacob and  S.T. Yau, {\em A special Lagrangian type equation for holomorphic line bundles}, Math. Ann. 369 (2017), no.1- 2, 869-898.

\bibitem{Jo} D. Joyce, {\em Conjectures on Bridgeland stability for Fukaya categories of Calabi-Yau manifolds, special Lagrangians, and Lagrangian mean curvature flow}, EMS Surv. Math. Sci. 2 (2015), no. 1, 1-62.


\bibitem{K} S. Ko\l{}odziej, {\em The complex Monge-Amp\`ere equation},  Acta Math. 180 (1998) 69--117.

\bibitem{LS} M. Lejmi and G. Szekelyhidi, {\em The J flow and stability},
Adv. Math. 274 (2015) 404-431.

\bibitem{LYZ}  N.C. Leung, S.T. Yau, and E. Zaslow, {\em From special Lagrangian to Hermitian-Yang-Mills via Fourier-Mukai}, Adv. Theor. Math. Phys. 4 (2000), no. 6, 1319-1341.

\bibitem{LT} P. Li and G. Tian, {\em On the heat kernel of the Bergmann metric on algebraic varieties}, J. Amer. Math. Soc. 8 (1995), no.4, 857-877.

\bibitem{Li} Y. Li, {\em Strominger-Yau-Zaslow conjecture for Calabi-Yau hypersurfaces in the Fermat family}, Acta Math. 229 (2022), no. 1, 1 - 53.

\bibitem{MMMS} M. Marino, R. Minasian, G. Moore, and A. Strominger,  {\em Nonlinear instantons from supersymmetric p-branes}, J. High Energy Phys. (2000), no. 1.

\bibitem{MeRa1} V.B. Mehta and A. Ramanathan, {\em Semistable sheaves over projective varieties and their restrictions to curves}, Math. Ann. 258 (1982) 213-224.

\bibitem{MeRa2} V.B. Mehta and A. Ramanathan, {\em Restriction of stable sheaves and representations of the fundamental group}, Inventiones Math. 77 (1984) 163-172.

\bibitem{Mi} M.L. Michelsohn, {\em On the existence of special metrics in complex geometry}, Acta Math. 149 (1982), 261-295.

\bibitem{PPZ1} D.H. Phong, S. Picard, and X.W. Zhang, {\em Geometric flows and the Strominger system}, Math. Z.  (2018) Vol. 288, 101-113.

\bibitem{PPZ2} D.H. Phong, S. Picard, and X.W. Zhang, {\em Anomaly flows}, Comm. Anal. Geom., Vol. 26, No. 4 (2018), 955-1008.


\bibitem{PPZ3} D.H. Phong, S. Picard, and X.W. Zhang,
{\em Fu-Yau equations}, J. Diff. Geom.  118 (2021) 147-187.

\bibitem{PhTo} D.H. Phong and D.T. To,
{\em Fully non-linear parabolic equations on compact Hermitian manifolds},
Ann. Sci. Ec. Normale Sup. 54 (3) (2021) 793-929, arXiv: 1711.10697.

\bibitem{Pi} V.P. Pingali, {\em The deformed Hermitian Yang-Mills equation on three-folds}, arXiv:1910.01870.


\bibitem{Po} D. Popovici, {\em A simple proof of a theorem by Uhlenbeck and Yau},
Math. Z. 250 (2005) 855-872.

\bibitem{Re} M. Reid, {\em The moduli space of $3$-folds with $K=0$ may nevertheless be irreducible}, Math. Ann. 278 (1987) no. 1-4, 329-334.

\bibitem{SaYu} O. Savin and H. Yu, {\em Regularity of optimal transport maps between planar convex domains}, Duke Math. J. 169 (2020) no. 7, 1305-1327.

\bibitem{ScYa} R. Schoen and S.T. Yau, {\em Lectures on differential geometry}, Conf. Proc. Lecture Notes Geom. Topology, I International Press, Cambridge, MA, 1994, v+235 pp.

\bibitem{Sh} B. Shiffman, {\em Complete characterization of holomorphic chains of codimension one}, Math. Ann. 274 (1986) 233-256.

 \bibitem{S20} J. Song, {\em Nakai-Moishezon criterions for complex Hessian equations}. arXiv:2012.07956.
 
\bibitem{SoWe} J. Song and B. Weinkove, {\em On the convergence and singularities of the 
$J$-flow with applications to the Mabuchi energy}, Commun. Pure Appl. Math. 61(2), 210-229 (2008).

\bibitem{St} A. Strominger, {\em Superstrings with torsion},
Nucl. Phys. B 274 (1986) no. 2, 253-284.

\bibitem{SYZ} A. Strominger, S.T. Yau, and E. Zaslow, {\em Mirror symmetry is T-duality}, Nuclear Phys. B 479 (1996), 243–259.

\bibitem{Stu} J. Sturm, private notes, available at 
{\texttt{https://sites.rutgers.edu/jacob-sturm/publications}}.

\bibitem{Sze} G. Sz\'ekelyhidi,  {\em Fully non-linear elliptic equations on compact hermitian manifolds}, J. Differential Geometry 109 (2018) 337-378.

\bibitem{ThYa} R.P. Thomas and  S.T. Yau, {\em Special Lagrangians, stable bundles and mean curvature flow} , Commun. Anal. Geom. 10 (2002), no. 5, 1075-1113.

\bibitem{Ti} G. Tian,  {\em On K\"ahler-Einstein metrics on certain K\"ahler manifolds with $C_1(M) > 0$}, Invent. Math. 89 (1987), no. 2, 225- 246.

\bibitem{TY90} G. Tian and S.T. Yau,  {\em Complete K\"ahler manifolds with zero Ricci curvature. I}, J. Amer. Math. Soc. 3 (1990), no. 3, 579 - 609.

\bibitem{TY89} G. Tian and  S.T. Yau,  {\em Complete K\"ahler manifolds with zero Ricci curvature, II}, Invent. Math., 106 (1): 27 - 60.

\bibitem{LiYa1} L.S. Tseng and S.T. Yau, {\em Cohomology and Hodge theory on symplectic manifolds: I}, J. Differential Geom. 91 (2012), no. 3, 383-416. 

\bibitem{LiYa2} L.S. Tseng and S.T. Yau, {\em Cohomology and Hodge theory on symplectic manifolds: II}, J. Differential Geom. 91 (2012), no. 3, 417-443.

\bibitem{LiYa3} L.S. Tseng and S.T. Yau, {\em Generalized cohomologies and supersymmetry}, Comm. Math. Phys. 326 (2014), no. 3, 875-885.

\bibitem{ToYa} F. Tong and S.T. Yau, {\em Generalized Monge-Amp\`ere functionals and related variational problems}, to appear in Amer. J. Math.,  arXiv:2306.01636.


\bibitem{UhYa} K. Uhlenbeck and S.T. Yau, S.T. {\em On the existence of Hermitian Yang-Mills-connections on stable bundles over K\"ahler manifolds}, Comm. Pure Appl. Math. 39 (1986) 257-293.

\bibitem{V} D. Vu, {\em Uniform diameter estimates for K\"ahler metrics}, arXiv:2405.14680.


\bibitem{Y} S.T. Yau, {\em On the Ricci curvature of a compact K\"ahler manifold and the complex Monge-Amp\`ere equation. I},  Comm. Pure Appl. Math. 31 (1978) 339-411.

\bibitem{Y80} S.T. Yau, {\em The role of partial differential equations in differential geometry}, Proceedings of the International Congress of Mathematicians (Helsinki, 1978),
237-250, Acad. Sci. Fennica, Helsinki, 1980.

}

\end{thebibliography}
\end{document}